\documentclass{article}

\usepackage{PRIMEarxiv}

\usepackage{cite}
\usepackage[utf8]{inputenc} 
\usepackage[T1]{fontenc}    
\usepackage{url}            
\usepackage{booktabs}       
\usepackage{amsfonts}       
\usepackage{amsmath}
\usepackage{nicefrac}       
\usepackage{microtype}      
\usepackage{lipsum}
\usepackage{fancyhdr}       
\usepackage{graphicx}
\usepackage{subcaption}
\usepackage[export]{adjustbox}
\usepackage{wrapfig}
\usepackage{algorithm}
\usepackage{algorithmic}
\usepackage{setspace}

\usepackage{tensor}  		

\newcommand{\mathtens}[1]{\tensor{\mathcal{#1}}{}}  	

\clearpage

\title{Tensor PCA from basis in tensor space}

\author{
  Claudio Turchetti and Laura Falaschetti\\
  DII - Department of Information Engineering \\
  Università Politecnica delle Marche\\
  via Brecce Bianche,12, 60131 Ancona, Italy\\
  \texttt{c.turchetti@univpm.it} \\
   \texttt{l.falaschetti@univpm.it} \\
}

\begin{document}

\setlength{\parindent}{20pt}

\large

\maketitle

\begin{abstract}
The aim of this paper is to present a mathematical framework for tensor PCA. The proposed approach is able to overcome the limitations of previous methods that extract a low dimensional subspace by iteratively solving an optimization problem.  The core of the proposed approach is the derivation of a basis in tensor space from a real self-adjoint tensor operator, thus reducing  the problem of deriving a basis to an eigenvalue problem. Three different cases have been studied to derive: i) a basis from a self-adjoint tensor operator; ii) a rank-1 basis; iii) a basis in a subspace. In particular, the equivalence between eigenvalue equation for a real self-adjoint tensor operator and standard matrix eigenvalue equation has been proven. For all the three cases considered, a subspace approach has been adopted to derive a tensor PCA.  Experiments on image datasets validate the proposed mathematical framework.
\end{abstract}



\section{Introduction}

\indent Tensors are arrays indexed by multiple indices, a generalization of vectors and matrices \cite{de1998matrix}, that have been adopted in several branches of data analysis, thanks  to their ability to represent a wide range of real-world data.
Tensors are most popular in machine learning 
\cite{ji2019survey}
\cite{tao2005supervised}, \cite{signoretto2014learning}, \cite{hou2017fast}, \cite{zafeiriou2009discriminant}, \cite{lebedev2014speeding}, \cite{kim2015compression}, \cite{lin2018holistic},  and image processing \cite{aja2009tensors}, \cite{sofuoglu2021multi}, \cite{qin2022low}, \cite{wang2021multi}, \cite{sun2022tensor},  \cite{long2021bayesian}, \cite{zhang2020robust}, \cite{tian2022low}, but they arise in numerous applications such as: sensor array processing \cite{sidiropoulos2000parallel}, big data representation \cite{kaur2018tensor}, completion of big data \cite{song2019tensor}, genomics \cite{omberg2007tensor}, neuroimaging data analysis \cite{zhou2013tensor}, \cite{sun2017store}, tumor classification 
\cite{sankaranarayanan2015tensor}, computer vision \cite{liu2012tensor}, finance data \cite{han2022rank}, human action and gesture recognition \cite{kim2007tensor}. 
Tensorial data are typically high dimensional involving high computational cost, thus techniques for dimensionality reduction are required to make possible the application of tensors to real data.
A widely adopted technique for this purpose is to represent high-dimensional tensors in a low-dimensional subspace, by preserving useful underlying information or structures. 
In this context several very powerful tensor decomposition approaches have been developed in the past \cite{kolda2009tensor}, \cite{cichocki2015tensor}, \cite{sidiropoulos2017tensor}, \cite{chen2021introduction}. Among these, Tucker decomposition \cite{tucker1963implications} and CANDECOMP/PARAFAC (CP) decomposition \cite{harshman1970foundations}, \cite{carroll1970analysis} are the most popular and fundamental models. Tucker model decomposes a tensor into a core tensor multiplied by a factor matrix along each mode, while CP model factorizes a tensor into a weighted sum of rank-1 tensors. 
Following these two seminal approaches, several extensions for low rank tensor decomposition have been proposed in the past few decades. Some of the most relevant developments in this field are: multilinear principal component analysis (MPCA), tensor rank-one decomposition (TROD), hierarchical Tucker PCA (HT-PCA) and tensor-train PCA (TT-PCA).\\
\indent
MPCA \cite{inoue2016generalized}, \cite{liu2010generalized}, \cite{lu2008mpca}, \cite{mavzgut2014dimensionality}, \cite{panagakis2009non}, \cite{sun2008incremental}, \cite{xu2008reconstruction}, \cite{yang2004two}, \cite{ye2004generalized}, \cite{ye2004gpca}, aims to find low-rank orthogonal projection matrices forming a lower dimensional subspace, such that the projection of a set of higher order data on this subspace minimizes a mean square error. The solution is iterative and performs feature extraction that captures most of the original tensorial input variations.\\
Similarly to MPCA, the TROD approach \cite{bro1997parafac}, \cite{faber2003recent}, \cite{harshman1970foundations}, \cite{kruskal1977three}, \cite{shashua2001linear} formulates the decomposition in terms of subspaces whose bases are rank-1 tensors. More formally, TROD aims to find a set of vectors forming a subspace, such that the projections of tensor data on this subspace minimizes a mean square error. An efficient greedy least square minimization procedure is used to iteratively compute each vector leaving all the others fixed.\\
\indent
Both HT and TT approaches fall within the class of tensor networks (TNs), that represents a higher order tensor as a set of sparsely interconnected lower order tensors, providing computational and storage benefits \cite{cichocki2014era}, \cite{cichocki2017tensor}, \cite{orus2014practical}.\\
\indent
HT-PCA \cite{grasedyck2010hierarchical},  \cite{grasedyck2013literature}, \cite{hackbusch2009new}, which has been proposed to reduce the memory requirements for Tucker decomposition, recursively splits the modes based on a hierarchy and creates a binary tree $T$ combining a subset of modes at each node. Thus, factor matrices are obtained from the singular value decomposition (SVD) of the matricization of a tensor corresponding to each node. In current applications  a suboptimal algorithm is used \cite{chaghazardi2017sample}, as estimating both tree and subspace is an NP hard problem.\\
\indent
TT-PCA decomposition \cite{bengua2017matrix}, \cite{chaghazardi2017sample}, is a special case of HT in which all nodes of the TN are connected in a cascade or train. Given a set of tensor data the objective of TT-PCA is to find order-3 tensors such that a left unfolding operator resulting in a matrix that represents a subspace is derived. In \cite{wang2019principal} an approach based on successive SVD algorithm for computing TT decomposition followed by thresholding the singular values, is proposed.\\
All the above methods for tensor decomposition try to find a low-dimensional subspace by solving an optimization problem. A commonly used approach for this purpose is the alternating least square (ALS) method \cite{zare2018extension}, which assumes the solution of all but one mode is known and then estimating the unknown parameter set of the remaining mode. However, even though this algorithm is efficient and easy to implement, it often converges to poor local minima and suffers from outliers and noise \cite{cheng2016probabilistic}.\\
\indent In order to overcome these issues, this paper proposes an alternative approach to tensor PCA that does not require solving an optimization problem.
Tensor PCA can be considered as multilinear extensions of principal component analysis (PCA)
\cite{jolliffe2016principal}, 
\cite{cunningham2015linear},
\cite{diamantaras1996principal},
\cite{shlens2014tutorial},
\cite{abdi2010principal},
\cite{sanguansat2012principal},
one of the most popular and efficient subspace learning techniques used in matrix analysis to reduce the dimensionality of a dataset. More formally, PCA derives an orthogonal linear transformation to extract a low-dimensional representation from high-dimensional data. This transformation can be regarded as a projection of data onto a coordinate system, termed as principal components (PCs), that generates the maximum variance.
There are two used approaches to perform the PCA: i) the eigenvalue decomposition of data covariance matrix; and ii) the SVD of data matrix. In both cases an orthogonal matrix $U$ is derived from the so-called Gram-matrix $G=X^TX$, where $X$ is the data matrix \cite{jolliffe2016principal}, while $X^T$ its transpose. The matrix $U$ represents the orthogonal linear transformation, whose columns are the eigenvectors of Gram-matrix and the existence of matrix $U$ is ensured by the symmetry of $G$. More specifically, the Gram-matrix can be interpreted as a real nonnegative self-adjoint operator that admits a set of real nonnegative eigenvalues, corresponding to a set of eigenvectors. 
Due to orthogonality of eigenvectors, the matrix $U$ represents a basis in the space of data vector $x$, whose observations are the rows of data matrix $X$. With these considerations in mind, the PCA can be contextualized in the more general framework of finite-dimensional linear spaces \cite{halmos2017finite} and in this context PCA can be considered as equivalent to derive a basis in vector space from a real self-adjoint operator. As this concept is quite general can be profitably used to extend PCA from vector-space to tensor-space. Thus in this paper some results of the theory of self-adjoint operators will be used as the foundation to derive tensor PCA.\\
\indent
The aim of this paper is to develop a mathematical framework for tensor PCA. The core of the proposed approach is the derivation of a basis in tensor space from a real \textit{self-adjoint operator}. It will be shown that tensors with a symmetrical structure belong to this important class of operators, can be easily derived as tensor products from tensor data, similar to the way in which Gram-matrix are derived from data matrix. In this way the problem of deriving a basis is equivalent to an eigenvalue problem. As final objective, tensors PCA is obtained as the projection of tensor data on the subspace spanned by the coordinates corresponding to the highest eigenvalues. In the paper three different problems will be carried out to derive: i) a basis from a real self-adjoint operator; ii) a rank-1 basis from a set of tensors; iii) a basis in a subset of tensors. For all these three cases a subspace approach will be adopted to derive a tensor PCA. In particular, with reference to the first issue, the properties of a self-adjoint tensor operator will be studied and the equivalence between eigenvalue equation for this operator and standard matrix eigenvalue equation will be proven. This result represents an advancement in the context of multilinear analysis, as it constitutes a generalization of previous results limited only to the case of supersymmetric tensor operators \cite{qi2005eigenvalues}, \cite{brachat2010symmetric}, \cite{cui2014all}, \cite{chang2009eigenvalue}, \cite{chen2022further}.\\
\indent
The paper is organized as follows.
Section 2 reports some preliminaries and notation. Section 3 deals with the finite linear space of tensors and some fundamental concepts to derive a tensor basis. In Section 4 the mathematical framework for tensor PCA is developed as projection of data on a subspace with reference to three different problems. 
Section 5 validates the proposed approach by experiments on image datasets.

\setlength{\parindent}{0pt}
\section{Preliminaries and notation}
\label{sec:preliminaries}

A real \textit{order-$d$ tensor} is a multidimensional array denoted by Euler script calligraphic letters, e.g. 
$\mathtens{X} \in \mathbb{R}^{I_1 \times I_2 \times \ldots \times I_d}$, where $\times$ represents the Cartesian product.
The number $d$ of dimensions, also known as \textit{modes}, is the order of a tensor.
The generic element of an order-$d$ tensor $\mathtens{X}$ will be denoted as

\begin{equation}\label{eq2.1}
\mathtens{X}_{\footnotesize\textbf{i}}=\mathtens{X}_{i_1, \ldots, i_d}, \;\; \textbf{i} = {(i_1, \ldots, i_d)}  \quad .
\end{equation}
or

\begin{equation}\label{eq2.2}
    \mathtens{X}(\textbf{i})=\mathtens{X}{(i_1, \ldots, i_d)} 
\end{equation}

where the index range in the $k$th mode is 
$1\leq{i_k}\leq{I_k},\hspace{0.2cm}  k=1:d$.
The vector $\textbf{i}$, in bold font, is a \textit{subscript vector} which assumes $L=I_1I_2\ldots{I_d}$ distinct values in the Cartesian interval $I=I_1 \times I_2 \times \ldots \times I_d$. $L$ represents the dimension of the $\mathtens{X}$ domain that will denoted by $dim(I)$.
\newline
\newline

\paragraph{Tensor product} Given two tensors

\begin{equation}\label{eq2.3}
    \mathtens{X}(\textbf{i})=x{(i_1, \ldots, i_d)}, \hspace{0.5cm} \mathtens{Y}(\textbf{j})=y{(j_1, \ldots, j_f)} 
\end{equation}

and assuming a common vector index $\textbf{k}$ exists such that  ${\textbf{i}}$, ${\textbf{j}}$ can be partitioned as 

\begin{equation}\label{eq2.4}
    \textbf{i}=(\textbf{l},\textbf{k},\textbf{m}), \hspace{0.5cm}\textbf{j}=(\textbf{p},\textbf{k},\textbf{q})
\end{equation}

the tensor product of $\mathtens{X}$ and $\mathtens{Y}$ along the multi-index $\textbf{k}$ combines the two tensors to give a third tensor

\begin{equation}\label{eq2.5}
\mathtens{Z}(\textbf{l},\textbf{m},\textbf{p},\textbf{q})=\sum_{\textbf{k=1}}^{\textbf{K}}\mathtens{X}(\textbf{l},\textbf{k},\textbf{m})\mathtens{Y}(\textbf{p},\textbf{k},\textbf{q})
\end{equation}

Here a multi-index notation is used for summation, meaning that, if $\textbf{k}={(k_1,\ldots,k_g)}$ is a \textit{length-g} index vector, then

\begin{equation}\label{eq2.6}
\sum_{\textbf{k=1}}^{\textbf{K}}=\sum_{{k_1}=1}^{{K_1}}...\sum_{{k_g}=1}^{{K_g}}
\end{equation}

As a result the product (\ref{eq2.6}) corresponds to a contraction of common index $\textbf{k}$.
Using the Einstein summation convention, that interprets repeated subscript as summation over that index, (\ref{eq2.6}) can be rewritten in a more compact form as

\begin{equation}\label{eq2.7}
\mathtens{Z}_{\footnotesize\hspace{0.02cm}\textbf{l},\textbf{m},\textbf{p},\textbf{q}}=
\mathtens{X}_{\footnotesize\hspace{0.02cm}\textbf{l},\textbf{k},\textbf{m}}
\mathtens{Y}_{\footnotesize\hspace{0.02cm}\textbf{p},\textbf{k},\textbf{q}}
\end{equation}

that implies summation over the vector index $\textbf{k}$. In the following we will extensively use the Einstein convention to simplify mathematical notation.
\newline
\paragraph{Inner product}  The inner product of two tensors of the same size $\mathtens{X}, \, \mathtens{Y} \in \mathbb{R}^{I_1 \times I_2 \times \ldots I_d}$ is a real scalar defined as

\begin{equation}\label{eq2.8}
\left\langle \mathtens{X}, \mathtens{Y} \right\rangle = \mathtens{X}_{\footnotesize\hspace{0.02cm}\textbf{i}}\mathtens{Y}_{\footnotesize\hspace{0.02cm}\textbf{i}}
\end{equation}
and $\mathtens{X},\mathtens{Y}$ are said to be orthogonal if $\left\langle \mathtens{X}, \mathtens{Y} \right\rangle=0$.
From this definition it follows that the norm of a tensor is given by 

\begin{equation}\label{eq2.9}
\|{\mathtens{X}}\| = \sqrt{\left\langle \mathtens{X}, \mathtens{X} \right\rangle}.
\end{equation}

\paragraph{Outer product} The outer product of tensor $\mathtens{X} = \mathtens{X}_{\footnotesize\textbf{i}} \in \mathbb{R}^{I_1 \times I_2 \times \ldots \times I_q}$ with
tensor $\mathtens{Y} = \mathtens{X}_{\footnotesize\textbf{j}} \in \mathbb{R}^{J_1 \times J_2 \times \ldots \times J_p}$ is the order-$(q+p)$ tensor $\mathtens{Z}$
defined as

\begin{equation}\label{eq2.10}
\mathtens{Z} = \mathtens{X} \circ \mathtens{Y}
\end{equation}

where the generic element of $\mathtens{Z}$ is 

\begin{equation}\label{eq2.11}
\mathtens{Z}_{\footnotesize\textbf{i},\textbf{j}} = 
\mathtens{X}_{\footnotesize\textbf{i}} \mathtens{Y}_{\footnotesize\textbf{j}},
\hspace{0.3cm} \textbf{i} = {(i_1, \ldots, i_q)},\hspace{0.2cm}\textbf{j}={(j_1, \ldots, j_p)}
\end{equation}

In particular for two vectors $x,y\in{R^{I_1}}$ the generic element of outer product ${Z} = x \circ y$ is the matrix

\begin{equation}\label{eq2.12}
z_{ij} = x_{i} \, y_{j} \quad .
\end{equation}

\paragraph{Vector-to-linear index transformation} This operation turns an order-$(p+q)$ tensor

\begin{equation}\label{eq2.13}
{\mathtens{X}}_{\footnotesize\textbf{j},\textbf{i}},\hspace{0.2cm}\textbf{j}
={(j_1, \ldots, j_p)},\hspace{0.2cm} \textbf{i} = {(i_1, \ldots, i_q)}
\end{equation}

to an order-$(p+1)$ tensor 

\begin{equation}\label{eq2.14}
{\mathtens{X}}_{\footnotesize\textbf{i},\large{m}},
\hspace{0.2cm}\textbf{j}={(j_1, \ldots, j_p)},
\hspace{0.2cm} m = 1:L 
\end{equation}

  transforming the vector index $\textbf{i}$ to the linear index $m$, ranging from 1 to $L$, (1:$L$ in Matlab notation) where $L=I_1I_2\ldots{I_q}$.

This is achieved by noticing that the indices $i_1, \ldots, i_q$ can be arranged in $L$ different ways (the number of permutations with repetition), so that a correspondence $\alpha$ between $m$ and $\textbf{i}$

\begin{equation}\label{eq2.15}
m=\alpha(\textbf{i}),\hspace{0.2cm}1\leq{m}\leq{L},
    \hspace{0.2cm}1\leq{i_k}\leq{I_k},
    \hspace{0.2cm}k=1:q
\end{equation}
holds.

The inverse transformation exists, since the correspondence is one-to-one, and is denoted by 

\begin{equation}\label{eq2.16}
    \textbf{i}=\alpha^{-1}(m),\hspace{0.2cm}1\leq{m}\leq{L},\hspace{0.5cm}1\leq{i_k}\leq{I_k},
    \hspace{0.2cm}k=1:q
\end{equation}

Among the possible choices for the correspondence, here we use the one described by the following relationship

\begin{equation}\label{eq2.17}
    m=\alpha(\textbf{i})=i_1+(i_2-1)I_1+(i_3-1)I_1I_2+\ldots+(i_q-1)I_1I_2\ldots{I_{q-1}}
\end{equation}

With reference to (\ref{eq2.17}) the matrix $T$ whose $m$th row represents the corresponding vector index $\textbf{i}(m)$ can be derived, so that the inverse transformation of (\ref{eq2.17}) for a value $m$ is formally given by

\begin{equation}\label{eq2.18}
    \textbf{i}=T(m,:).
\end{equation}

An example of such a matrix for $\textbf{i}=(i_1,i_2,i_3)$,\hspace{0.1cm}$I_1=3,I_2=2,I_3=2$ is the following

\begin{center}
$T=
\left( 
\begin{array}{c c c c c c c c c c c c}
 1 & 1 & 1 & 1 & 1 & 1 & 2 & 2 & 2 & 2 & 2 & 2 \\ 
 1 & 1 & 1 & 2 & 2 & 2 & 1 & 1 & 1 & 2 & 2 & 2 \\ 
 1 & 2 & 3 & 1 & 2 & 3 & 1 & 2 & 3 & 1 & 2 & 3 \\ 
\end{array}
\right)^T $
\end{center}

where $1\leq{m}\leq{12}$ represents the row-index.

\section{The space of tensors}
\label{sec:space of tensors}

It is straightforward to show that the set of order-$d$ tensors $\mathtens{X} \in \mathbb{R}^{I_1 \times I_2 \times \ldots \times I_d}$ forms a linear space, since it is closed under addition and multiplication by a scalar.
With reference to the canonical basis $e_{i_1} = (1, 0, \ldots, 0), \ldots, e_{i_d} = (0,0,\ldots,1)$ in vector space,
any order-$d$ tensor $\mathtens{X}$ can be decomposed as

\begin{equation}\label{eq3.1}
\mathtens{X} =  \mathtens{X}_{i_1 i_2 \ldots i_d} \; e_{i_1} e_{i_2} \ldots e_{i_d}
\end{equation}

The outer product $e_{i_1} e_{i_2} \ldots e_{i_d}$ is equivalent to a set of $L = \prod_{j=1}^{d} I_j$ order-$d$ tensors, thus due to property (\ref{eq3.1}) it represents a canonical tensor basis.
As this basis is of size $L$, the set of order-$d$ tensors form a linear space of \textit{dimensionality} $L$.
As an example for $d=2$  and $I_1 = I_2 = 3$ we have $L=9$  and

$ e_1 e_1 =
\left( 
\begin{array}{ccc}
1 & 0 & 0 \\
0 & 0 & 0 \\
0 & 0 & 0 
\end{array} 
\right), \; \nonumber $

$e_1 e_2 =
\left( 
\begin{array}{ccc}
0 & 1 & 0 \\
0 & 0 & 0 \\
0 & 0 & 0 
\end{array} 
\right), \; \nonumber $
$e_1 e_3 =
\left( 
\begin{array}{ccc}
0 & 0 & 1 \\
0 & 0 & 0 \\
0 & 0 & 0 
\end{array} 
\right), \; \ldots$


\paragraph{Tensor operators} In the product (\ref{eq2.7}) the tensor $\mathtens{X}_{\footnotesize\hspace{0.02cm}\textbf{l},\textbf{k},\textbf{m}}$ can be interpreted as a transformation from space $R^{P\times{K}\times{Q}}$ to space $R^{L\times{M}\times{P}\times{Q}}$, where $P=p_1\times\ldots\times{p_p}$, $K=k_1\times\ldots\times{k_k}$ and so on,
transforming tensor
$\mathtens{Y}_{\footnotesize\hspace{0.02cm}\textbf{p},\textbf{k},\textbf{q}}$ 
to tensor $\mathtens{Z}_{\footnotesize\hspace{0.02cm}\textbf{l},\textbf{m},\textbf{p},\textbf{q}}$. Among the possible transformations, we will interested in transformations from
$\mathbb{R}^{I_1 \times I_2 \times \ldots \times I_d}$ to itself. In this particular case we have 

\begin{equation}\label{eq3.2}
{\mathtens{Z}}_{\footnotesize\hspace{0.02cm}\textbf{i}}={\mathtens{A}}_{\footnotesize\hspace{0.02cm}\textbf{i},\textbf{j}}{\mathtens{Y}}_{\footnotesize\hspace{0.02cm}\textbf{j}}
\end{equation}

where
$\textbf{i} = {(i_1, \ldots, i_d)}$ and $\textbf{j} = {(j_1, \ldots, j_d)}$, 
so that the tensors  $\mathtens{Z}_{\footnotesize\hspace{0.02cm}\textbf{i}}$ and $\mathtens{Y}_{\footnotesize\hspace{0.02cm}\textbf{j}}$ are both in $\mathbb{R}^{I_1 \times I_2 \times \ldots \times I_d}$ and  the order-$2d$ tensor $\mathtens{A}_{\footnotesize\hspace{0.02cm}\textbf{i},\textbf{j}}$ is said to be a tensor \textit{operator}. Thus, denoting with $\mathtens{L}_d$ the order-$d$ tensor space $\mathbb{R}^{I_1 \times I_2 \times \ldots \times I_d}$ the operator $\mathtens{A}_{\textbf{i,j}}$  establishes a transformation from
$\mathtens{L}_d$ to $\mathtens{L}_d$ and (\ref{eq3.2}) can be rewritten without ambiguity as

\begin{equation}\label{eq3.3}
\mathtens{Z} = \mathtens{A} \mathtens{Y}
\end{equation}

A tensor $\mathtens{V}$ in $\mathtens{L}_d$ is said to be an \textit{eigentensor} if $\mathtens{V}\neq{0}$ and if for some scalar $\lambda$ the following equation

\begin{equation}\label{eq3.4}
\mathtens{A}\mathtens{V}=\lambda\mathtens{V}
\end{equation}
is satisfied.
The scalar $\lambda$ is known as the \textit{eigenvalue} of $\mathtens{A}$ associated with the eigentensor $\mathtens{V}$,

\paragraph{Real self-adjoint operators}
\label{sec:theory} 
An important class of operators is the class of real \textit{self-adjoint} operators $\mathtens{A}$, for which the following property

\begin{equation}\label{eq3.7}
    \left\langle \mathtens{A}\mathtens{Y}, \mathtens{Z} \right\rangle =\left\langle \mathtens{Y}, \mathtens{A}\mathtens{Z} \right\rangle 
\end{equation}
holds. 

An order-$2d$ tensor 
\begin{equation}\label{eq3.5}
{\mathtens{A}}_{\footnotesize\hspace{0.02cm}\textbf{i},\textbf{j}},
\hspace{0.2cm}\textbf{i}={(i_1, \ldots, i_d)},\hspace{0.2cm} \textbf{j} = {(j_1, \ldots, j_d)}
\end{equation}

such that

\begin{equation}\label{eq3.6}
{\mathtens{A}}_{\footnotesize\hspace{0.02cm}\textbf{i},\textbf{j}}={\mathtens{A}}_{\footnotesize\hspace{0.02cm}\textbf{j},\textbf{i}}
\end{equation}

is a real self-adjoint operator, since it results 

\begin{equation}\label{eq3.8}
    \left\langle \mathtens{A}\mathtens{Y}, \mathtens{Z} \right\rangle ={\mathtens{A}}_{\footnotesize\hspace{0.02cm}\textbf{i},\textbf{j}}{\mathtens{Y}}_{\footnotesize\hspace{0.02cm}\textbf{j}}{\mathtens{Z}}_{\footnotesize\hspace{0.02cm}\textbf{i}}={\mathtens{Y}}_{\footnotesize\hspace{0.02cm}\textbf{j}}{\mathtens{A}}_{\footnotesize\hspace{0.02cm}\textbf{j},\textbf{i}}{\mathtens{Z}}_{\footnotesize\hspace{0.02cm}\textbf{i}}=\left\langle \mathtens{Y}, \mathtens{A}\mathtens{Z} \right\rangle.
\end{equation}

It is straightforward to show that the covariance of tensor $\mathtens{X}_{\textbf{i}}$ defined as 

\begin{equation}
\mathtens{R}_{\textbf{i,j}}=E\{
\mathtens{X}_{\textbf{i}}
\mathtens{X}_{\textbf{j}}\},
\end{equation}
where $E\{.\}$ denotes the mean value, is a self-adjoint operator.
In the particular case of order-1 tensor space , achieved with $d=1$ in (\ref{eq3.5}), the self-adjoint operator (\ref{eq3.6}) reduces to a \textit{symmetric} matrix that operates on space $R^{I_1}$. \\
On the basis of this property it follows that eigentensors of a self-adjoint operator $\mathtens{A}$ belonging to distinct eigenvalues are orthogonal.
In fact suppose that $\mathtens{A}\mathtens{V}_1=\lambda_1\mathtens{V}_1$ and $\mathtens{A}\mathtens{V}_2=\lambda_2\mathtens{V}_2$  for $\lambda_1\neq\lambda_2$, if $\mathtens{A}$ is self-adjoint , then

\begin{equation}\label{eq3.9}
    \left\langle \mathtens{A}\mathtens{V}_1, \mathtens{V}_2 \right\rangle =\lambda_1\left\langle \mathtens{V}_1,\mathtens{V}_2 \right\rangle
\end{equation}
and also

\begin{equation}\label{eq3.10}
    \left\langle \mathtens{A}\mathtens{V}_1, \mathtens{V}_2 \right\rangle =\left\langle \mathtens{V}_1, \mathtens{A}\mathtens{V}_2\right\rangle =\lambda_2\left\langle \mathtens{V}_1,\mathtens{V}_2 \right\rangle
\end{equation}

Therefore,  $(\lambda_1-\lambda_2)\left\langle \mathtens{V}_1,\mathtens{V}_2 \right\rangle=0$, and hence $\left\langle \mathtens{V}_1,\mathtens{V}_2 \right\rangle=0$, since $\lambda_1\neq\lambda_2$.

The eigenvalues of a real self-adjoint operator are real, indeed combining (\ref{eq3.4}) and (\ref{eq3.7}) $\left\langle \mathtens{V}, \mathtens{A}\mathtens{V} \right\rangle =
\lambda{\left\langle \mathtens{V}, \mathtens{V} \right\rangle}=\lambda$
that yields a real scalar

\begin{equation}\label{eq3.11}
    \lambda=\left\langle \mathtens{V}, \mathtens{A}\mathtens{V} \right\rangle/
    {\left\langle \mathtens{V}, \mathtens{V} \right\rangle}
\end{equation}

From (\ref{eq3.11}) it follows that operators for which $\left\langle \mathtens{V}, \mathtens{A}\mathtens{V} \right\rangle\geq0 $ holds, ensure the eigenvalues are nonnegative. Operators belonging to a such class are said to be \textit{nonnegative operators}. A nonnegative operator can be easily derived from a tensor $\mathtens{X}_{\textbf{i,j}}$ as follows

\begin{equation}\label{eq3.12}
\mathtens{A}_{\textbf{i,j}}=
\mathtens{X}_{\textbf{k,i}}
\mathtens{X}_{\textbf{k,j}}
\end{equation}

For such an operator we have

\begin{equation}\label{eq3.13}
\left\langle \mathtens{V}, \mathtens{A}\mathtens{V} \right\rangle=
\mathtens{V}_{\textbf{i}}
\mathtens{A}_{\textbf{i,j}}
\mathtens{V}_{\textbf{j}}=
\mathtens{X}_{\textbf{k,i}}
\mathtens{V}_{\textbf{i}}
(\mathtens{X}_{\textbf{k,j}}
\mathtens{V}_{\textbf{j}})=
\mathtens{Z}_{\textbf{k}}
\mathtens{Z}_{\textbf{k}}\geq0
\end{equation}
where $\mathtens{Z}_{\textbf{k}}=
\mathtens{X}_{\textbf{k,i}}
\mathtens{V}_{\textbf{i}}$,
so that from (\ref{eq3.11}) it results $\lambda\geq0$. The operator $\mathtens{A'}_{\textbf{i,j}}=
\mathtens{X}_{\textbf{i,k}}
\mathtens{X}_{\textbf{j,k}}$ satisfies this property as well.\\ 
In particular for linear indices $\textbf{i}=m$,\hspace{0.2cm}$\textbf{j}=n$,\hspace{0.2cm}$\textbf{k}=l$, the operator reduces to a matrix 
\begin{equation*}
    {G}_{m,n}={X}_{l,m}{X}_{l,n}=
    ({X}^T_{m,l}){X}_{l,n}
\end{equation*}
that is $G=X^TX$. Similarly, operator 
$\mathtens{A'}_{\textbf{i,j}}$ reduces to the matrix $G^T=XX^T$. The matrix $G$ represents the well known Gram-matrix, while matrix $G^T$ its transpose.
Thus, tensor $\mathtens{A}_{\textbf{i,j}}$ is the analog in tensor space of Gram-matrix $G$.\\
Tensor eigenvalues and eigenvectors have received much attention recently in the literature \cite{qi2007eigenvalues},
\cite{qi2012spectral},
\cite{chang2013survey}, following the seminal works of Qi \cite{qi2005eigenvalues}, and Lim \cite{lim2005singular}, that introduced spectral theory of \textit{supersymmetric} tensors. A  tensor $\mathtens{A}$ is called supersymmetric as its entries $a_{i_1,\dots{i_m}}$ are invariant under any permutations of their indices.
It is worth to note that tensor (\ref{eq3.6}) has a symmetrical structure but it is not supersymmeric, thus in order to distinguish it from the class of supersymmetric tensors, simply referred as \textit{"symmetric"} in some papers \cite{brachat2010symmetric}, \cite{cui2014all}, \cite{chang2009eigenvalue}, \cite{chen2022further}, in the following we will adopt the term self-adjoint for tensor 
(\ref{eq3.6}).
    

\section{Tensor PCA from basis in tensor space}
\label{sec:tensor PCA}

As discussed in the Introduction the core aspect of principal component analysis (PCA) is the derivation of a basis from a self-adjoint operator, from which a representation on a subspace that minimizes an error condition can be achieved.
The aim of this Section is to derive a \textit{basis} in the order-$d$ tensor space $\mathbb{R}^{I_1 \times I_2 \times \ldots \times I_d}$, that is a set of $L=I_1I_2\ldots{I_d}$ orthogonal tensors

\begin{equation}\label{eq4.1}
\mathtens{U}_{\footnotesize\textbf{i},l},
\hspace{0.2cm}l=1:L
\end{equation}

such that, for any tensor $\mathtens{X}_{\footnotesize\textbf{i}}$ in this space, the following orthogonal expansion

\begin{equation}\label{eq4.1.1}
\mathtens{X}_{\footnotesize\textbf{i}}=
 \mathtens{U}_{\footnotesize\textbf{i},l}d_l
\end{equation}

holds. Using the orthogonality property   
${\mathtens{U}_{\footnotesize\textbf{i},l}}
{\mathtens{U}_{\footnotesize\textbf{i},l'}}=
\delta_{l,l'}$ of basis ($\delta_{l,l'}$ is the Dirac operator: $\delta_{l,l'}=1$ for $l=l'$, while $\delta_{l,l'}=0$ for $l\neq{l'}$ ), the vector of coefficients $d_l$ in (\ref{eq4.1.1}) can be easily derived by the inner product of $\mathtens{X}_{\footnotesize\textbf{i}}$ and 
${\mathtens{U}_{\footnotesize\textbf{i},l}}$

\begin{equation}\label{eq4.1.2}
d_l=\mathtens{X}_{\footnotesize\textbf{i}}
 \mathtens{U}_{\footnotesize\textbf{i},l}
\end{equation}

As final step tensor PCA is derived by extracting low-dimensional subspace representation of dataset that ensures a minimum reconstruction error.
In the following  results  will be presented with reference to three different problems to derive:
i) a basis assuming a real self-adjoint tensor operator is given; ii) a rank-1 basis from a set of tensors; iii) a basis in a subset of tensors.


\subsection*{A. Basis from real self-adjoint operator}
\label{sec:theory_basis}
Let us given the order-$2d$ tensor $\mathtens{A}$

\begin{equation}\label{eq4A1}
    \mathtens{A}_{\footnotesize\textbf{i,j}},
    \hspace{0.2cm} \textbf{i} = {(i_1, \ldots, i_d)},
    \hspace{0.2cm}\textbf{j}={(j_1, \ldots, j_d)}
\end{equation}
    
with $1\leq{i_k}\leq{I_k}, \hspace{0.2cm}  1\leq{j_k}\leq{J_k},\hspace{0.2cm}k=1:d$, so that it results

\begin{equation}\label{eq4A2}
    \mathtens{A} \in \mathbb{R}^{I \times J},\hspace{0.2cm}{I=I_1 \times I_2 \times \ldots \times I_d},\hspace{0.2cm} {J=J_1 \times J_2 \times \ldots \times J_d}
\end{equation}

and $dim(I)=dim(J)=I_1I_2\ldots{I_d}=L$. Assuming $\mathtens{A}$ is self-adjoint, i.e. $\mathtens{A}_{\footnotesize\textbf{i,j}}=
\mathtens{A}_{\footnotesize\textbf{j,i}}$,
then we are interested in solving the following eigenvalue problem, that is to determine the tensors $\mathtens{U}$ (the eigentensors) and the scalars $\lambda$ (the eigenvalues) such that equation 

\begin{equation}\label{eq4A3}  
\mathtens{A}\mathtens{U}=\lambda\mathtens{U}
\end{equation}

is satisfied, where $\mathtens{U} \in \mathbb{R}^{I}$. To this end we will prove the following proposition.


\paragraph{Proposition1} Eq. (\ref{eq4A3}) is equivalent to the eigenvalue equation

\begin{equation}\label{eq4A4}
  a(n,m)u(m)=\lambda{u(n)},\hspace{0.2cm}m,n=1:L
\end{equation}

where $a(m,n)$ is a symmetric matrix achieved using the vector-to-linear index transformation (\ref{eq2.17}) and the eigentensors $\mathtens{U}$ are obtained by applying the inverse vector-to-linear transformation  (\ref{eq2.18}) to the eigenvectors $u$ in (\ref{eq4A4}). \\

\paragraph{Proof} Using index convention introduced in Section 2,  Eq. (\ref{eq4A3}) is equivalent to the $L$ equations

\begin{equation}\label{eq4A5}   
\mathtens{A}_{\footnotesize\textbf{i,j}}
 \mathtens{U}_{\footnotesize\textbf{j}}=
 \lambda\mathtens{U}_{\textbf{i}},
 \hspace{0.3cm}\textbf{i}\in{I},
 \hspace{0.2cm}L=dim(I)
\end{equation}

which can be written as 

\begin{equation}\label{eq4A6}
(\mathtens{A}_{\footnotesize\textbf{i,j}}  -\lambda\delta_{\footnotesize\textbf{i,j}})
\mathtens{U}_{\footnotesize\textbf{j}}=0,
\end{equation}
where the term on the left 
$\mathtens{Y}_{\textbf{i}}=
(\mathtens{A}_{\footnotesize\textbf{i,j}}  -\lambda\delta_{\footnotesize\textbf{i,j}})
\mathtens{U}_{\footnotesize\textbf{j}}
$ represents an order-d tensor. The solutions of (\ref{eq4A6}) are independent on how the terms in the equation are arranged as entries in the tensor $\mathtens{Y}_{\textbf{i}}$, meaning that (\ref{eq4A6}) is invariant to a one-to-one transformation of subscript vectors.
Assuming indices $\textbf{i},\textbf{j}$ are obtained by the inverse vector-to-linear index transformation (\ref{eq2.18}) then, substituting the vector indices
$\textbf{i}(n)=T(n,:)$,
$\textbf{j}(m)=T(m,:)$, (\ref{eq4A6})  becomes

\begin{equation}\label{eq4A7}
(\mathtens{A}_{\footnotesize\textbf{i}(n),
    \footnotesize\textbf{j}(m)}
    -\lambda\delta_{\footnotesize\textbf{i}(n),\textbf{j}(m)})\mathtens{U}_{\footnotesize\textbf{j}(m)}=0
\end{equation}

By defining the $(L\times{L})$ matrices

\begin{equation}\label{eq4A8}
a(n,m)=\mathtens{A}_{\footnotesize\textbf{i}(n),\textbf{j}(m)}, 
    \hspace{0.2cm} \delta(n,m)=\delta_{\footnotesize\textbf{i}(n),\textbf{j}(m)}
\end{equation}

and the vector

\begin{equation}\label{eq4A9} 
u(m)=\mathtens{U}_{\footnotesize\textbf{j}(m)},
\end{equation}

thus (\ref{eq4A7}) can be re-arranged as 

\begin{equation}\label{eq4A10}
    (a(n,m))-\lambda\delta(n,m))u(m)=0
\end{equation}

or equivalently

\begin{equation}\label{eq4A11}
 a(n,m)u(m)=\lambda{u(n)},\hspace{0.2cm}m,n=1:L
\end{equation}
As $\mathtens{A}_{\textbf{i,j}}$ is self-adjoint, then it results 

\begin{equation}\label{eq4A12}
    a(n,m)=\mathtens{A}_{\footnotesize\textbf{i}(n),\textbf{j}(m)}=
    \mathtens{A}_{\footnotesize\textbf{j}(m),\textbf{i}(n)}=
    a(m,n)
\end{equation}

that proves matrix $a(n,m)$ is symmetric (or self-adjoint), and this concludes the proof. \\
Proposition1 shows that by a one-to-one vector-to-linear transformation the eigenvalue tensor equation (\ref{eq4A3}) reduces to an eigenvalue matrix equation.
(\ref{eq4A4}) is a standard eigenvalue equation for matrix $a(n,m)$, that admits $L$ orthogonal eigenvectors due to the symmetry of $a(n,m)$.  Once the eigenvectors $u(n)$ are obtained by solving (\ref{eq4A4}), they can be converted to tensors $\mathtens{U}$ of proper size using the inverse vector-to-linear index transformation (\ref{eq2.17})

\begin{equation}\label{eq4A13}
 u(n(\textbf{i))}=
 \mathtens{U(\textbf{i})}=
 \mathtens{U}_{\footnotesize\textbf{i}}
\end{equation}.

The statement proven with Proposition 1 shows that an eigenvalue equation for real self-adjoint tensor operators is invariant under one-to-one transformations of subscript vectors, thus reducing to a matrix eigenvalue equation. This result represents an advancement in the general higher order tensor eigenvalue problem, so far limited to the class of supersymmetric tensor operators \cite{qi2005eigenvalues}, \cite{lim2005singular}\cite{qi2007eigenvalues},\cite{qi2012spectral},
\cite{chang2013survey}.

Having derived a basis in the order-$d$ tensor space, a generic tensor $\mathtens{X}_{\footnotesize\textbf{i}} \in \mathbb{R}^{I_1 \times I_2 \times \ldots \times I_d}$ in this space can be represented by the orthogonal expansion (\ref{eq4.1.1}). A pseudo-code of the procedure previously described to derive a basis from a self-adjoint operator is described in 
Algorithm~\ref{alg:symmetricOperator}.

\paragraph{ Tensor PCA}
On the basis of previous results a tensor PCA can be carried out by deriving a low-dimensional subspace representation of dataset that minimizes the reconstruction error. The following proposition establishes this result.

\paragraph{Proposition 2} - Let us refer to a generic basis 
$\{{\mathtens{U}_
{\footnotesize\textbf{i},l},\hspace{0.2cm}
l=1:L}\}$ such that for any
$\mathtens{X}_{\footnotesize\textbf{i}} \in \mathbb{R}^{I_1 \times I_2 \times \ldots \times I_d}$ the representation

\begin{equation*}
\mathtens{X}_{\footnotesize\textbf{i}}=
 \mathtens{U}_{\footnotesize\textbf{i},l}d_l
\end{equation*}
\begin{equation}\label{eq4.15}
d_l=\mathtens{X}_{\footnotesize\textbf{i}}
 \mathtens{U}_{\footnotesize\textbf{i},l}
\end{equation}
holds, and decompose $\mathtens{X}_{\footnotesize\textbf{i}}$ as

\begin{equation}\label{eq4A14}
\mathtens{X}_{\footnotesize\textbf{i}}=
{\hat{\mathtens{X}}_{\footnotesize\textbf{i}}+\mathtens{N}_{\mathbf{i}}}
\end{equation}

where $\hat{\mathtens{X}}_{\footnotesize\textbf{i}}=\mathtens{U}_{\footnotesize\textbf{i},l}d_l,\hspace{0.2cm}l=1:M$ with $M<L$, represents a subspace representation and $\mathtens{N}_{\mathbf{i}}=
\mathtens{U}_{\footnotesize\textbf{i},p}d_p,\hspace{0.2cm}p=M+1:L$ the residual.
Then the basis $\{{\mathtens{U}_
{\footnotesize\textbf{i},l},\hspace{0.2cm}
l=1:L}\}$ that minimizes the error

\begin{equation}\label{eq4A15}
\mathtens{E}=E{\{\|\mathtens{X}_
{\textbf{i}}-
{\hat{\mathtens{X}}_{i}}\|^2}\}
\end{equation}

is solution of the eigenvalue equation 

\begin{equation}\label{eq4A16}
\mathtens{R}_{\textbf{i,j}}
\mathtens{U}_
{\footnotesize\textbf{j},p}=\lambda_p
{\mathtens{U}_
{\footnotesize\textbf{i},p}}
\end{equation}

where $\mathtens{R}_{\textbf{i,j}}=
E\{\mathtens{X}_
{\textbf{i}}\mathtens{X}_{\textbf{j}}\}$ 
 is the covariance of 
$\mathtens{X}_{\textbf{i}}$.\\

\vspace{0.2cm}
A proof of this proposition is given as follows.

\paragraph{Proof}  From $\mathtens{X}_{\footnotesize\textbf{i}}=
{\hat{\mathtens{X}}_{\footnotesize\textbf{i}}+\mathtens{N}_{\mathbf{i}}}$ the mean square error $\mathtens{E}$ results as given by

\begin{equation*}
\mathtens{E}=E{\{\mathtens{N}_
{\textbf{i}}
{{\mathtens{N}}_{\textbf{i}}}}\}=
E{\{\mathtens{U}_
{\textbf{i},p}d_{p}
{{\mathtens{U}}_{\textbf{i},q}}d_{q}\}}=
\end{equation*}
\begin{equation*}
E{\{\delta_{p,q}d_{p}d_{q}\}}=
E{\{\mathtens{X}_{\textbf{i}}\mathtens{U}_
{\textbf{i},p}
{{\mathtens{X}_{\textbf{j}}\mathtens{U}}_
{\textbf{j},q}}\}}=
\end{equation*}
\begin{equation}
\mathtens{U}_{\textbf{i},p}
E\{\mathtens{X}_{\textbf{i}}
\mathtens{X}_{\textbf{j}}\}
\mathtens{U}_{\textbf{j},q}=
\mathtens{U}_{\textbf{i},p}
\mathtens{R}_{\textbf{i,j}}
\mathtens{U}_{\textbf{j},p},
\end{equation}
with $p=M+1:L$.
Thus to find the basis that minimizes the error $\mathtens{E}$ subject to the orthonormal condition for the basis $\mathtens{U}_{\textbf{i},p}$, is equivalent to the problem

\begin{equation}
\min \{\mathtens{U}_{\textbf{i},p}
\mathtens{R}_{\textbf{i,j}}
\mathtens{U}_{\textbf{j},p}\},
\hspace{0.2cm}
s.t. 
\hspace{0.2cm}
\mathtens{U}_{\textbf{i},p}
\mathtens{U}_{\textbf{i},p}=1
\end{equation}
Using the Lagrange multiplier method the problem reduces to the equation

\begin{equation}
\frac{\partial}{\partial{\mathtens{U}_{\textbf{i},p}}}
\{\mathtens{U}_{\textbf{i},p}
\mathtens{R}_{\textbf{i,j}}
\mathtens{U}_{\textbf{j},p}+
\lambda_p
{(1-\mathtens{U}_{\textbf{i},p}
\mathtens{U}_{\textbf{i},p})}
\}=0
\end{equation}

which once solved yields
\begin{equation}
\mathtens{R}_{\textbf{i,j}}
\mathtens{U}_{\textbf{j},p}=
\lambda_p\mathtens{U}_{\textbf{i},p}
\end{equation}
that is the eigenvalue equation of covariance $\mathtens{R}_{\textbf{i,j}}$, and this concludes the proof.\\
As a corollary of this proposition, it follows that the basis elements in the truncated representation $\hat{\mathtens{X}}_{\footnotesize\textbf{i}}=\mathtens{U}_{\footnotesize\textbf{i},l}d_l$, $l=1:M$ correspond to the highest eigenvalues in (\ref{eq4A16}).
In fact the error

\begin{equation}
\min{\mathtens{E}}=
\mathtens{U}_{\textbf{i},p}
(\mathtens{R}_{\textbf{i,j}}
\mathtens{U}_{\textbf{j},p})=
\mathtens{U}_{\textbf{i},p}
\lambda_p\mathtens{U}_{\textbf{i},p}=
\end{equation}
\begin{equation}
\lambda_{p}\|\mathtens{U}_{\textbf{i},p}\|^2=
\sum_{p=M+1}^{L}\lambda_p
\end{equation}

is minimum when the eigenvalues corresponding to the basis elements in residual are the lowest or, that is the same, when the eigenvalues corresponding to the basis elements in the subspace are the highest. \\
Proposition 2 proves that low-dimensional subspace representation $\hat{\mathtens{X}_{\textbf{i}}}$ is optimal, since minimizes the mean square reconstruction error.

\begin{algorithm}[!htb]
\large
\caption{Basis from real self-adjoint operator}
\label{alg:symmetricOperator}
\setstretch{1.1} 
\begin{algorithmic}

\STATE{INPUT: order-$2d$ self-adoint operator}
\STATE{$\mathtens{A}_{\footnotesize\textbf{i,j}},
\hspace{0.2cm} \textbf{i} = {(i_1, \ldots, i_d)},
\hspace{0.2cm}\textbf{j}={(j_1, \ldots, j_d)}$}
\vspace{0.2cm}

\STATE{1. Derive ($L\times{d}$) matrix $T$ defining vector-to-linear transformation}
\STATE{$\textbf{i}(n)=T(n,:),\hspace{0.2cm}
\textbf{j}(m)=T(m,:),\hspace{0.2cm}m,n=1:L$}
\vspace{0.2cm}

\STATE{2. Compute matrix}
\STATE{$a(n,m)=\mathtens{A}_{\footnotesize\textbf{i}(n),\textbf{j}(m)},\hspace{0.2cm}m,n=1:L$}
\vspace{0.2cm}

\STATE{3. Solve eigenvalue equation}
\STATE{$a(n,m)u(m)=\lambda{u(n)},\hspace{0.2cm}m,n=1:L$}
\vspace{0.2cm}

\STATE{4. Convert eigenvector $u(n)$ to tensor by inverse vector-to-linear transformation}
\STATE {$\mathtens{U}_{\footnotesize\textbf{i}}=
 u(n(\textbf{i}))$}
\vspace{0.2cm}

\STATE{OUTPUT: the basis}
\STATE{$\mathtens{U}=\{\mathtens{U}_{\mathbf{i},l},\hspace{0.2cm}l=1:L\}$}
\end{algorithmic}
\end{algorithm}


\subsection*{B. Rank-1 basis}
In this case, a set 

\begin{equation}\label{eq4B1}
    \mathtens{X}=\{\mathtens{X}_{\footnotesize\textbf{i},n},\hspace{0.2cm}n=1:N\}
\end{equation}

of $N$ order-$d$ tensors
$\mathtens{X}_{\footnotesize\textbf{i}}\in\mathbb{R}^{I_1 \times I_2 \times \ldots \times I_d},$
$\hspace{0.2cm} \textbf{i} = (i_1, \ldots, i_d)$, is given. The problem to be solved is to derive a set of 
$L=I_1I_2\ldots{I_d}$, rank-$1$ orthogonal tensors 
$\mathtens{U}_{\footnotesize\textbf{i},\footnotesize\textbf{j}}$ 
with ${\textbf{j} = (j_1, \ldots, j_d)}$
having the same dimension of {\textbf{i}}, which can be written as

\begin{equation}\label{eq4B2}
    \mathtens{U}_{\footnotesize\textbf{i,j}}=
    {U}_{i_{1},j_{1}}
    {U}_{i_{2},j_{2}}\ldots{U}_{i_{d},j_{d}}
\end{equation}

where the generic element 
${U}_{i_{k},j_{k}}\in\mathbb{R}^{I_k \times J_k}$
is an orthonormal matrix $U^k$ that corresponds to the $k$-mode. In such a way for a given value of $\textbf{j}$ the resulting tensor

\begin{equation}\label{eq4B3}
    {U}_{\footnotesize\textbf{:j}}=
    {U}_{:,j_{1}}
   {U}_{:,j_{2}}\ldots{U}_{:,j_{d}}
\end{equation}

is a rank-$1$ tensor being expressed as outer product of $d$ vectors, i.e. the columns of matrices 
${U}_{i_{k}j_{k}}$.
The problem can be addressed by deriving an eigenvector equation for each mode of the tensor $\mathtens{X_{\textbf{i}}}$, so that once solved an orthogonal matrix is obtained. This is equivalent to find a set of self-adjoint tensors that operates along the $d$ modes. To this end let us refer to the $k$th mode. The tensor 

\begin{equation}\label{eq4B4}
    \mathtens{X}(i_1,\ldots,{i_k},\ldots,{i_d},n)
    \mathtens{X}(i_1,\ldots,{j_k},\ldots,{i_d},n)=A^k(i_k,j_k)
\end{equation}

achieved by the contraction of all indices but $i_k$ and $j_k$, defines a self-adjoint operator along $k$th mode as a symmetric matrix $A^k$ of dimension $(I_k\times{J_k})$.  From this matrix we can derive an orthogonal matrix $U^k$ which is solution of the eigenvector equation

\begin{equation}\label{eq4B5}
    A^kU^k=U^k\Sigma^k
\end{equation}

whose generic element is ${U}_{i_{k}j_{k}}$. By repeating this approach for all the modes of tensor, we finally obtain the required $d$ orthogonal matrices that can be used in (\ref{eq4B2}) to derive  a set of L orthogonal rank-$1$ tensors 

\begin{equation}\label{eq4B6}
    \mathtens{U}_{\footnotesize\textbf{i,j}},
    \hspace{0.2cm} \textbf{i} = (i_1, \ldots, i_d),
    \hspace{0.2cm} \textbf{j} = (j_1, \ldots, j_d)
\end{equation}

Transforming vector index $\textbf{i}$ to linear index by (\ref{eq2.18}) yields the set

\begin{equation}\label{eq4B7}
\hat{\mathtens{{U}}}_{\footnotesize\textbf{i},m}=
\mathtens{U}_{\footnotesize\textbf{i},\footnotesize\textbf{j}{(m)}}, \hspace{0.2cm}, m=1:L
\end{equation}

that is a basis for the space of tensors $\mathtens{X}_{\footnotesize\textbf{i}}\in\mathbb{R}^{I_1 \times I_2 \times \ldots \times I_d}$. As a consequence each tensor of data set $\mathtens{X}$ can be represented as

\begin{equation}\label{eq4B8}
\mathtens{X}_{\footnotesize\textbf{i},n}=
\hat{\mathtens{U}}_{\footnotesize\textbf{i},\large{m}}d_{n,m}, \hspace{0.2cm}n=1:N,\hspace{0.2cm} m=1:L
\end{equation}

where the $(N\times{L})$ matrix ${d_{n,m}=[D]_{n,m}}$ is given by the inner product

\begin{equation}\label{eq4B9}
d_{n,m}=\mathtens{X}_{\footnotesize\textbf{i},n}
\hat{\mathtens{U}}_{\footnotesize\textbf{i},{m}}
\end{equation}

The matrix $D$ can be decomposed by the SVD as

\begin{equation}\label{eq4B10}
D=Y\Sigma{Z}^T    
\end{equation}
where

\begin{equation}\label{eq4B11}
 \Sigma =
\left( 
\begin{array}{cc}
S & 0  \\
0 & 0  
\end{array} 
\right), \; 
\hspace{0.2cm}S=diag(\sigma_1,\ldots,\sigma_r)
\end{equation}
and $r=rank(D)$.
The generic term of $D$ is

\begin{equation}\label{eq4B12}
d_{n,m}=Y_{n,l}S_{l,p}Z_{m,p},
\hspace{0.2cm}l,p=1:r,
\end{equation}
then, from (\ref{eq4B8})

\begin{equation}\label{eq4B13}
\mathtens{X}_{\footnotesize\textbf{i},n}=
\hat{\mathtens{U}}_{\footnotesize\textbf{i},m}
Y_{n,l}S_{l,p}Z_{m,p} 
\end{equation}

and $S_{l,p}=\delta_{l,p}\sigma_p$, it results

\begin{equation}\label{eq4B14}
\mathtens{X}_{\footnotesize\textbf{i},n}=
\hat{\mathtens{U}}_{\footnotesize\textbf{i},m}
\sigma_{l}Y_{n,l}Z_{m,l},
\hspace{0.2cm}l=1:r
\end{equation}

where $\sigma_{l}Y_{n,l}Z_{m,l}$ is the rank-1 decomposition of matrix $D$, and $\hat{\mathtens{U}}_{\footnotesize\textbf{i},m}\in{R^{I\times{L}}}$ is the rank-1 basis.

A pseudo-code of the procedure previously described to derive a rank-1 basis is described in 
Algorithm~\ref{alg:rank-1 basis}.

\begin{algorithm}[!htb]
\large
\caption{ Rank-1 basis}
\label{alg:rank-1 basis}
\setstretch{1.1} 
\begin{algorithmic}

\STATE{INPUT: dataset of tensors}
\STATE{$\mathtens{X}=\{\mathtens{X}_{\footnotesize\textbf{i},n},\hspace{0.2cm}n=1:N\},\hspace{0.2cm} \textbf{i} = (i_1, \ldots, i_d),
\hspace{0.2cm}1\leq{i_k}\leq{I_k}$}
\vspace{0.2cm}
\STATE{for k=1:d}

\STATE{1. Derive the $(I_k\times{I_k})$ symmetric matrix $A^k$ from tensor product along $k$-mode.}
\STATE{$A^k(i_k,j_k)=\mathtens{X}
(i_1,\ldots,{i_k},\ldots,{i_d},n)
\mathtens{X}(i_1,\ldots,{j_k},\ldots,{i_d},n)$}
\vspace{0.2cm}

\STATE{2. Compute matrix $U^k=U_{i_k,j_k}$ solving eigenvalue equation}
\STATE{$A^k=U^k\Sigma^kU^k$}
\STATE{end}
\vspace{0.2cm}

\STATE{3.The basis as outer product of matrices $U_{i_k,j_k}$}
\STATE{$\mathtens{U}_{\footnotesize\textbf{i,j}}=
{U}_{i_{1},j_{1}}
{U}_{i_{2},j_{2}}
\ldots{U}_{i_{d},j_{d}}$}
\vspace{0.2cm}

\STATE{4. Transform the index $\textbf{i}$ to linear index $m$ }
\STATE{$\hat{\mathtens{{U}}}_{\footnotesize\textbf{i},m}=
\mathtens{U}_{\footnotesize\textbf{i},\footnotesize\textbf{j}{(m)}}, \hspace{0.2cm}, m=1:L$}
\vspace{0.2cm}

\STATE{5. Derive the $(N\times{L})$ matrix $(D)_{n,m}=d_{n,m}$ as inner product}
\STATE{$d_{n,m}=\mathtens{X}_
{\footnotesize\textbf{i},n}
\hat{\mathtens{U}}_{\footnotesize\textbf{i},m}$}
\vspace{0.2cm}

\STATE{6. Decompose the matrix $d_{n,m}$ by SVD}
\STATE{$d_{n,m}=Y_{n,l}S_{l,p}Z_{m,p}
\hspace{0.2cm}l,p=1:rank(D)$}
\vspace{0.2cm}

\STATE{OUTPUT: orthogonal expansion of tensor $\mathtens{X}_{\footnotesize\textbf{i},n}$}
\STATE{$\mathtens{X}_
{\footnotesize\textbf{i},n}=
\hat{\mathtens{U}}_
{\footnotesize\textbf{i},\large{m}}
\sigma_{l}Y_{n,l}Z_{m,l}, 
\hspace{0.2cm} l=1:r$}
\end{algorithmic}
\end{algorithm}


\paragraph{Tensor PCA} A principal component analysis can be carried out as follows. As the eigenvalues are nonnegative, the elements of $S_{l,p}$ are put in decreasing order and the corresponding eigenvectors in $Y_{n,l}$ and $Z_{m,p}$ are ordered accordingly. Then the  principal components analysis corresponds to the truncated representation 

\begin{equation}\label{eq4B15}
\hat{\mathtens{X}}_
{\footnotesize\textbf{i},n}=
\hat{\mathtens{U}}_
{\footnotesize\textbf{i},\large{l}}
d_{n,l}, 
\hspace{0.2cm} l=1:M
\end{equation}

with $M<r$, that minimizes the error

\begin{equation}\label{eq4B16}
\mathtens{E}=\|\mathtens{X}_
{\textbf{i},n}-
{\hat{\mathtens{X}}_{\textbf{i},n}\|^2}=
\hat{\mathtens{U}}_
{\footnotesize\textbf{i},\large{m}}d_{n,m}
\hat{\mathtens{U}}_
{\footnotesize\textbf{i},\large{m}}d_{n,m}=
d_{n,m}d_{n,m}
\end{equation}
with $m=M+1:r$. Using (\ref{eq4B12}) and orthonormality of $Y_{n,l}$ and $Z_{m,l}$ we have 

\begin{equation}\label{eq4B17}
\mathtens{E}=Y_{n,l}Y_{n,l'}
\sigma_{l}\sigma_{l'}
Z_{m,l}Z_{m,l'}=\sigma_{l}\sigma_{l},
\hspace{0.2cm}l=M+1:r.
\end{equation}
The error is minimum when the singular values $\sigma_l$ corresponding to the  basis elements in residual are the lowest, or when the singular values corresponding to elements in truncated representation are the highest. Thus (\ref{eq4B15}) is an optimal subspace representation of dataset $\mathtens{X}$.

\subsection*{C. Basis of a subspace}

Let us suppose a set 

\begin{equation}\label{eq4C1}
 \mathtens{X}=\{\mathtens{X}_{\footnotesize\textbf{i},n},
 \hspace{0.2cm}n=1:N\}
\end{equation}

of N order-$d$ tensors $\mathtens{X}_{\footnotesize\textbf{i}},$
is given, with 
$N<L=I_1I_2\ldots{I_d}$. Assuming $r$ tensors within this set are linearly independent, thus this  subset spans a subspace in tensor space 
$\mathbb{R}^{I_1 \times I_2 \times \ldots \times I_d}$, 
whose dimension is $r\leq{N}$. The problem of deriving a basis from this subset is equivalent to the QR problem \cite{golub2013matrix}, that is well known in matrix algebra, which aims to derive an orthonormal set of vectors, from a set of independent vectors. In tensor space the problem can be formulated as follows. Assume every tensors in $\mathtens{X}$ can be represented as linear combination of $r$ orthogonal tensors
$\mathtens{Q}_{\textbf{i},l},\hspace{0.1cm}l=1:r$. Formally we have

\begin{equation}\label{eq4C2}
\mathtens{X}_{\footnotesize\textbf{i},n}=
\mathtens{Q}_{\textbf{i},l}r_{n,l},
\hspace{0.2cm}l=1:r
\end{equation}
where $r_{n,l}$ is the $(N\times{r})$ matrix of coefficients in the linear combination (\ref{eq4C2}).

The tensor product of $\mathtens{X}_
{\footnotesize\textbf{i},n}$ with itself along the vector index $\textbf{i}$, gives 
\newline

\begin{equation}\label{eq4C5}
\mathtens{X}_{\footnotesize\textbf{i},n}
\mathtens{X}_{\footnotesize\textbf{i},m}  =\mathtens{Q}_{\textbf{i},l}r_{n,l}
 \mathtens{Q}_{\textbf{i},l'}r_{m,l'}=
 r_{n,l}r_{m,l'}\delta_{l,l'}=
r_{n,l}r_{m,l},\hspace{0.2cm}
\end{equation}

where the orthogonality of 
$\mathtens{Q}_{\textbf{i},l}$, i.e.
$\mathtens{Q}_{\textbf{i},l}
\mathtens{Q}_{\textbf{i},l'}=
\delta_{l,l'}$, is used. The tensor product (\ref{eq4C5}) represents a real nonnegative self-adjoint operator in the space $R^N$, due to symmetry of the $(N\times{N})$ matrix $r_{n,l}r_{m,l}$, thus ensuring a set of nonnegative eigenvalues exists.

It can be easily shown that the rank of matrix $[R]_{n,l}=r_{n,l}$ is $r$. In fact due to the orthogonality of $\mathtens{Q}_{\textbf{i},l}$, 
it results 
\newline

\begin{equation}\label{eq4C3}
    \mathtens{Q}_{\textbf{i},l'}
    \mathtens{X}_{\footnotesize\textbf{i},n}=
    \mathtens{Q}_{\textbf{i},l'}
    \mathtens{Q}_{\textbf{i},l}r_{n,l}=
    \delta_{l,l'}r_{n,l}=r_{n,l}
\end{equation}
and 

\begin{equation*}
r_{n,l}r_{n,l'}=\mathtens{Q}_{\textbf{i},l}
\mathtens{X}_{\footnotesize\textbf{i},n}
\mathtens{Q}_{\textbf{i},l'}
\mathtens{X}_{\footnotesize\textbf{i},n}=
\end{equation*}
\begin{equation}\label{eq4C4}
     \big\|\mathtens{X}\big\|^2
     \mathtens{Q}_{\textbf{i},l}
    \mathtens{Q}_{\textbf{i},l'}=
    \big\|\mathtens{X}\big\|^2
    \delta_{l,l'}
\end{equation}

(\ref{eq4C4}) shows that two generic columns $r_{n,l}$ and 
$r_{n,l'}$ of matrix $R$ are orthogonal each other and, as a consequence $r=rank(R)$.\\
Due to symmetry, the $(N\times{N})$ matrix $r_{n,l}r_{m,l}$ can be decomposed as

\begin{equation*}
    r_{n,l}r_{m,l}=
    U_{n,p}\Sigma^2_{p,q}U^T_{q,m}=
    U_{n,p}\Sigma^2_{p,q}U_{m,q}=
\end{equation*}
\begin{equation}\label{eq4C6}
 U_{n,r}\Sigma_{p,q}\Sigma_{p,q}U_{m,q}=
    U_{n,p}\Sigma_{p,l}(U_{m,q}\Sigma_{q,l})
\end{equation}

The matrix on the left $r_{n,l}r_{m,l}=[RR^T]_{n,m}$ is the Gram-matrix of $R^T$, and it is well known from matrix algebra that $rank(RR^T)=rank(R)$. This implies that the matrix $\Sigma_{p,q}$ is an $(r'\times{r'})$ diagonal matrix $\Sigma=diag(\sigma_1\ldots{\sigma_r})$, whose elements $\sigma_j$ are the singular values, with $r'=rank(R^TR)=rank(R)=r$. It is worth to notice that any other matrix $RV$, whith $V$ orthogonal, satisfies (\ref{eq4C6}), thus $RR^T$ can be decomposed in a non unique way as 

\begin{equation}\label{eq4C7}
    r_{n,l}=U_{n,p}\Sigma_{p,l}
\end{equation}

Combining (\ref{eq4C2}) and (\ref{eq4C7}) gives
\begin{equation}\label{eq4C8}
    \mathtens{X}_{\footnotesize\textbf{i},n}=
    \mathtens{Q}_{\textbf{i},l}
    U_{n,p}\Sigma_{p,l}
\end{equation}
and the tensor product of $\mathtens{X}$ with matrix $U$ yields

\begin{equation*}
 \mathtens{X}_{\footnotesize\textbf{i},n}U_{n,m}=
 \mathtens{Q}_{\textbf{i},l}\Sigma_{p,l}
 U_{n,p}U_{n,m}=
 \end{equation*}
 \begin{equation}\label{eq4C9}
 \mathtens{Q}_{\textbf{i},l}\Sigma_{p,l}
 \delta_{p,m}=
 \mathtens{Q}_{\textbf{i},l}\Sigma_{l,m}
\end{equation}

due to the orthogonality of $U$.
The matrix $\Sigma_{p,l}$ is diagonal full rank, thus $\Sigma^{-1}_{p,l}$ exists and from (\ref{eq4C9})
\begin{equation}\label{eq4C10}
    \mathtens{X}_{\footnotesize\textbf{i},n}
    U_{n,m}\Sigma^{-1}_{m,l}=
    \mathtens{Q}_{\textbf{i},l}
\end{equation}
or
\begin{equation}\label{eq4C11}
    \mathtens{Q}_{\textbf{i},l}=
    \mathtens{X}_{\footnotesize\textbf{i},n}
    b_{n,l},\hspace{0.2cm}l=1:r
\end{equation}
where
\begin{equation}\label{eq4C12}
    b_{n,l}=U_{n,m}\Sigma^{-1}_{m,l}
\end{equation}
is an $(N\times{N})$ matrix. (\ref{eq4C11}) and (\ref{eq4C12}), together with the decomposition of $RR^T$ give the requested relationships for the basis of subspace $\mathtens{X}$.

A pseudo-code of the procedure previously described to derive a basis of a subspace is described in 
Algorithm~\ref{alg:subspaceBasis}.

\begin{algorithm}[!htb]
\large
\caption{Basis of subspace}
\label{alg:subspaceBasis}
\setstretch{1.1} 
\begin{algorithmic}

\STATE{INPUT: dataset of tensors}
\STATE{$\mathtens{X}=\{\mathtens{X}_{\footnotesize\textbf{i},n},
 \hspace{0.2cm}n=1:N\}$}
\vspace{0.2cm}
\STATE{for k=1:d}

\STATE{1. Derive the symmetric $(N\times{N})$ matrix $(A)_{m,n}=a_{m,n}$ as tensor product }
\STATE{$a_{m,n}=\mathtens{X}_
{\footnotesize\textbf{i},m}
\mathtens{X}_{\footnotesize\textbf{i},n}$}
\vspace{0.2cm}

\STATE{2. Decompose the matrix $a_{m,n}$ by solving eigenvalue equation}
\STATE{$a_{m,n}=U_{n,p}\Sigma^2_{p,q}U_{m,q}
\hspace{0.2cm}p,q=1:rank(A)$}
\vspace{0.2cm}

\STATE{3. Compute the matrix $b_{n,l}$}
\STATE{$ b_{n,l}=U_{n,m}\Sigma^{-1}_{m,l}$}
\vspace{0.2cm}

\STATE{OUTPUT: The basis}
\STATE{$\mathtens{Q}_{\textbf{i},l}=
    \mathtens{X}_{\footnotesize\textbf{i},n}
    b_{n,l},\hspace{0.2cm}l=1:rank(A)$}
\end{algorithmic}
\end{algorithm}

\paragraph{Tensor PCA} To derive a principal component analysis, as the eigenvalues are nonnegative, we put the elements of $\Sigma_{p,l}$  in decreasing order and we order the corresponding eigenvectors in $U_{n,p}$  accordingly. Then by defining the truncated matrices

\begin{equation}\label{eq4C13}
    \hat{U}_{n,p}=U_{n,p},\hspace{0.2cm}
\hat{\Sigma}_{n,p}=\Sigma_{n,p},\hspace{0.2cm}
p,l=1:M
\end{equation}
with $M<r$, principal components analysis corresponds to the subspace representation 

\begin{equation}\label{eq4C14}
 \hat{\mathtens{X}}_{\textbf{i},n}=
{{Q}}_{\textbf{i},l}
\hat{U}_{n,p}\hat{\Sigma}_{p,l}\hspace{0.2cm}
n=1:N,\hspace{0.2cm}p,l=1:M
\end{equation}
that minimizes the error

\begin{equation}\label{eq4C15}
\mathtens{E}=\|\mathtens{X}_
{\textbf{i},n}-
{\hat{\mathtens{X}}_{\textbf{i},n}\|^2}=
{{Q}}_{\textbf{i},l}
\hat{U}_{n,p}\hat{\Sigma}_{p,l}
{{Q}}_{\textbf{i},l'}
\hat{U}_{n,p'}\hat{\Sigma}_{p',l'}
\end{equation}
Due to orthonormality of ${{Q}}_{\textbf{i},l}$ and $\hat{U}_{n,p}$ it results

\begin{equation}\label{eq4C16}
\mathtens{E}=
\hat{\Sigma}_{p,l}\hat{\Sigma}_{p,l}=\sigma_p\sigma_p
\end{equation}

with  $p=M+1:r$. Thus the error is minimum when the singular values corresponding to the  basis elements in the residual are the lowest, or when the singular values corresponding to elements in truncated representation are the highest. Thus (\ref{eq4C14}) is an optimal subspace repesentation of database $\mathtens{X}$.

\section{Experimental results }
\label{sec:experiment}

The aim of the experiments is to validate the three approaches previously described in Section 4 for tensor PCA. To this end we  refer to the space of order-3 tensors, in which the generic element represents an RGB image. The experiments have been conducted on three different datasets, namely TinyImageNet, CALTECH-101 and FLOWERS, to show that tensor PCA gives a low-dimensional subspace representation of datasets. 

\subsection{Experiments on TinyImageNet}

The first experiment aims to validate the approach discussed in Subsection 4.A, in which a basis is derived from a self-adjoint operator.
In this experiment the TinyImageNet dataset has been used\cite{Le2015TinyIV, ChrabaszczLH17Tiny, tiny-imagenet}. 
The dataset contains 100000 images of 200 classes (500 for each class) downsized to 64×64 colored images, which are a subset of the ImageNet dataset \cite{Deng2009ImageNet}. Each class has 500 training images, 50 validation images and 50 test images.
The experiment refers to the  reconstruction of the 100000 RGB images of the  dataset by representation (\ref{eq4.1.1}). 
For this purpose, the data has been organized in a order-$4$ tensor $\mathtens{X} = 
\{ \mathtens{X}_{\textbf{i},n},
n=1:N\}$, with
$N = 100000, \; \mathtens{X} \in \mathbb{R}^{64 \times 64 \times 3\times100000}$, so that the dimension of the basis is $L = 64 \cdot 64 \cdot 3 = 12288$.
A self-adjoint operator can be derived from the dataset $\mathtens{X}$ as tensor product

\begin{equation*}
    \mathtens{A}_{\textbf{i},\textbf{j}}=
    \mathtens{X}_{\textbf{i},n}
    \mathtens{X}_{\textbf{j},n}
\end{equation*}

since it results in $\mathtens{A}_{\textbf{i},\textbf{j}}=\mathtens{A}_{\textbf{j},\textbf{i}}$, thus satisfying the requirements for Proposition 1. Solving equation (\ref{eq4A5}) by Algorithm 1 a set of $L$ eigenvalues and eigentensors are achieved.

\begin{figure}[ht]
\includegraphics[width=0.9\linewidth,height=8cm]
{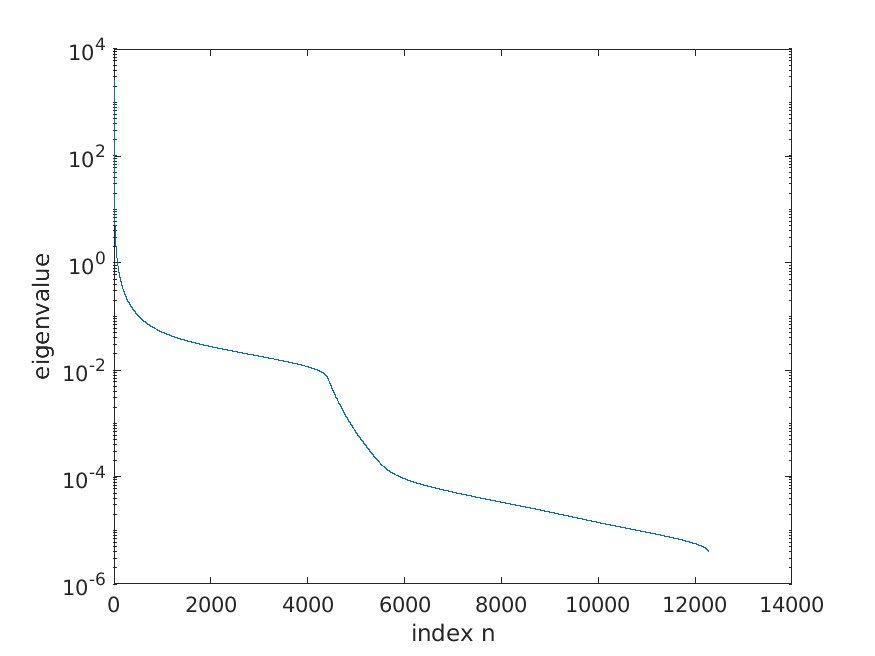}
\caption{\large{The behavior of eigenvalues as obtained from TinyImageNet dataset.}}
\label{fig:TINYeigenvalues}
\end{figure}

\begin{figure}[ht]
\includegraphics[width=0.9\linewidth,height=8cm]
{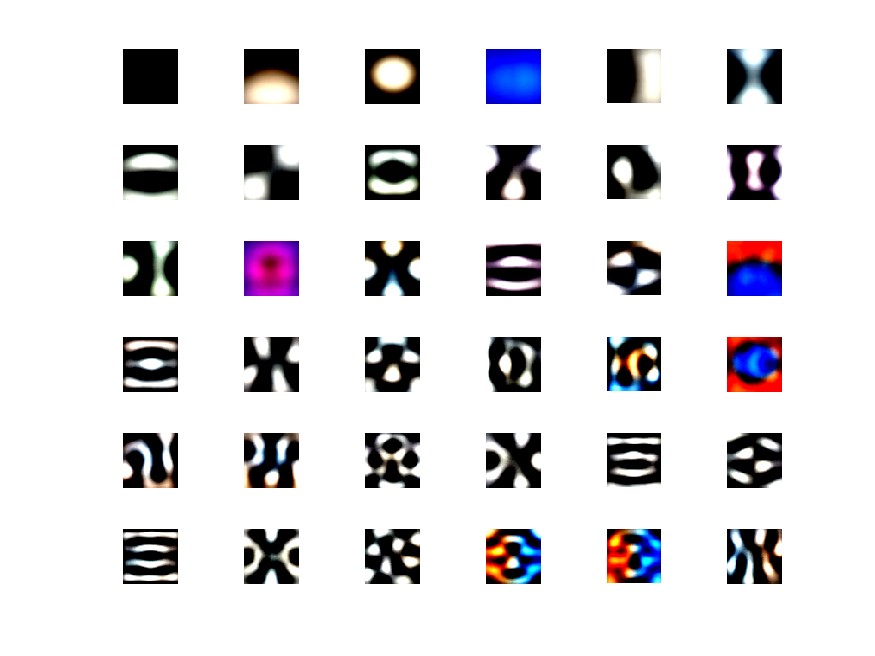}
\caption{\large{Some basis elements from TinyImageNet dataset corresponding to the highest eigenvalues.}}
\label{fig:TINYbasis}
\end{figure}

\begin{figure}[ht]
\includegraphics[width=0.9\linewidth,height=8cm]
{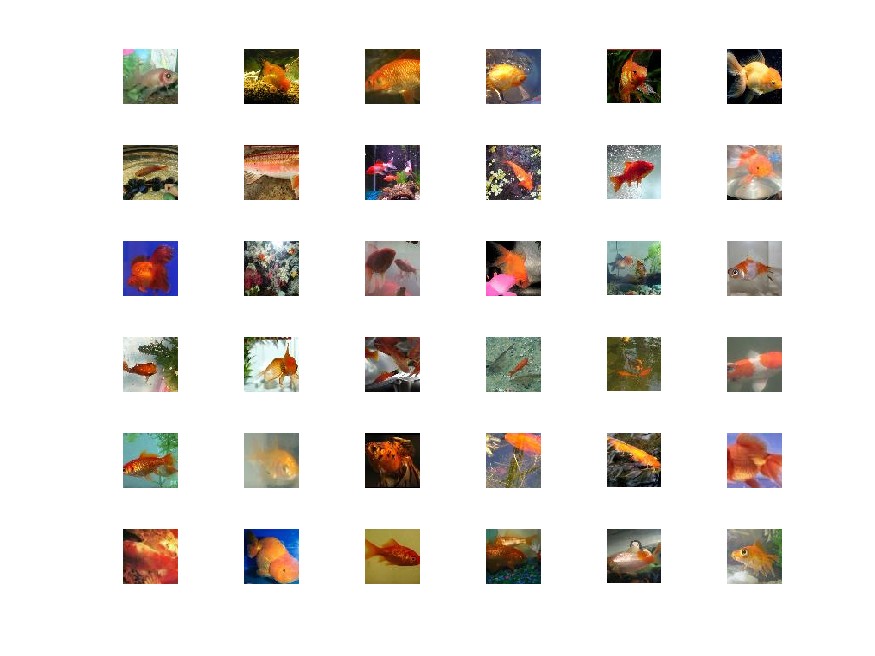}
\caption{\large{Several images taken from the dataset TinyImageNet.}}
\label{fig:TINY}
\end{figure}

\begin{figure}[ht]
\includegraphics[width=0.9\linewidth,height=8cm]
{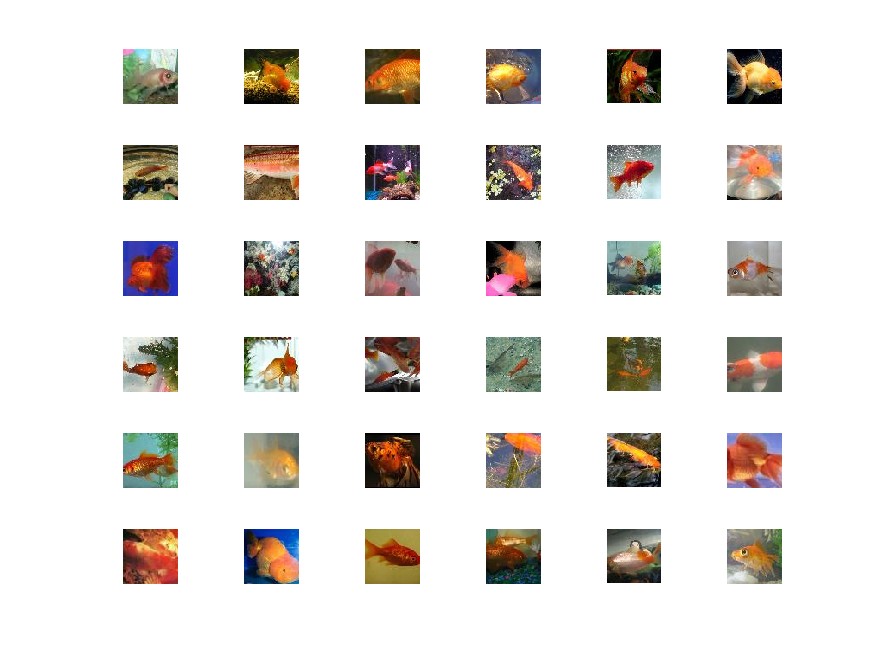}
\caption{\large{The same images of Fig.~\ref{fig:TINY} reconstructed by (\ref{eq4.1.1}).}}
\label{fig:TINYreconstructed}
\end{figure}

\begin{figure}[ht]
\includegraphics[width=0.9\linewidth,height=8cm]{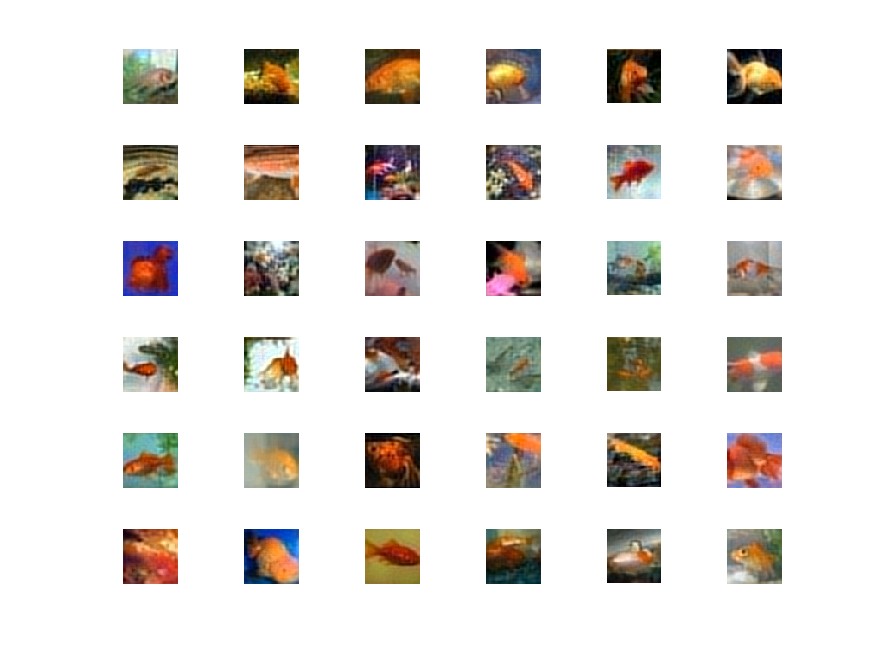}
\caption{\large{The same images of Fig.~\ref{fig:TINY} reconstructed with truncated basis.}}
\label{fig:TINYtruncatedRecon}
\end{figure}

Fig.~\ref{fig:TINYeigenvalues} shows the behavior of eigenvalues as obtained from TinyImageNet dataset, using Algorithm 1,
while some basis elements from TinyImageNet dataset corresponding to the first eigenvalues are reported in Fig.~\ref{fig:TINYbasis}.
Fig.~\ref{fig:TINY}
reports several images taken from the dataset. For comparison the same images reconstructed by (\ref{eq4.1.1}), with the basis achieved with Algorithm 1, are reported in Fig.~\ref{fig:TINYreconstructed}.
As you can see the two datasets match perfectly, thus proving the validity of the approach.
On the basis of the behavior of eigenvalues depicted in Fig.~\ref{fig:TINYeigenvalues}
tensor PCA has been applied retaining only the principal components, i.e. the ones corresponding to the most significant eigenvalues. 
Fig.~\ref{fig:TINYtruncatedRecon} shows the same images as Fig.~\ref{fig:TINY}, reconstructed with the 500 most significant components of the basis, thus confirming the correctness of approximation by tensor PCA.

\subsection{Experiments on CALTECH-101}

The second experiment refers to the Algorithm 2 discussed in Subsection 4.B, in which a rank-$1$ basis is derived from a dataset.
For this purpose CALTECH-101 dataset 
\cite{caltech101}, that contains images from 101 object categories (e.g., “helicopter”, “elephant” and “chair” etc.) and a background category that contains the images not from the 101 object categories, has been used.  For each object category, there are about 40 to 800 images, while most classes have about 50 images. The resolution of the image is roughly about 300×200 pixels.
In this experiment 2990 RGB images reduced to a single size
of ${64 \times 64 \times 3}$ have been used for image reconstruction by the representation (\ref{eq4.1.1}). For this purpose, the data has been organized in an order-$4$ tensor $\mathtens{X} = 
\{ \mathtens{X}_{\textbf{i},n},n=1:N\}$, with $N = 2990, \; \mathtens{X} \in \mathbb{R}^{64 \times 64 \times 3\times2990}$, so that the dimension of the basis is $L = 64 \cdot 64 \cdot 3 = 12288$.

\begin{figure}[ht]
\includegraphics[width=0.9\linewidth,height=8cm]
{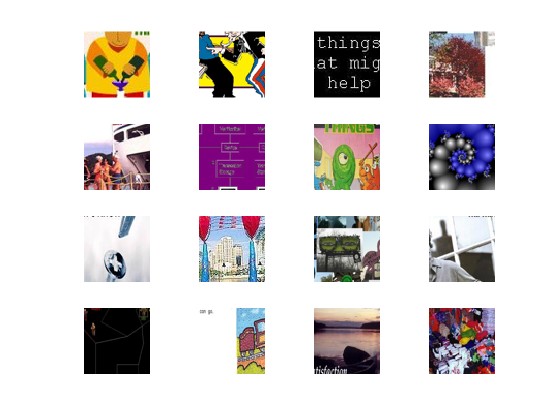}
\caption{\large{Several images taken from the dataset CALTECH-101.}}
\label{fig:caltech101}
\end{figure}

\begin{figure}[ht]
\includegraphics[width=0.9\linewidth,height=7cm]{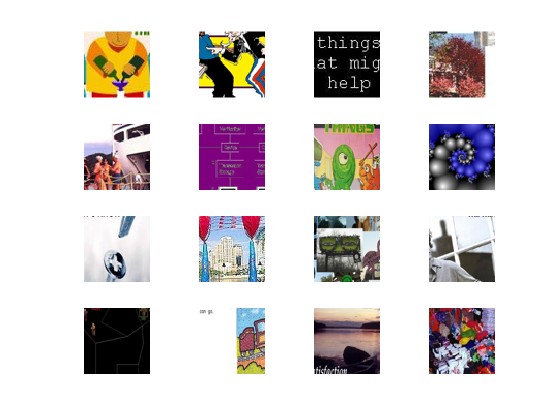}
\caption{\large{The same images of Fig.~\ref{fig:caltech101} reconstructed by (\ref{eq4.1.1}).}}
\label{fig:caltech101reconstructed}
\end{figure}

Fig.~\ref{fig:caltech101}
reports several images taken from the dataset. For comparison the same images reconstructed by (\ref{eq4.1.1}), with the basis achieved with Algorithm 2, are reported in Fig.~\ref{fig:caltech101reconstructed}.
Also in this case these results show a perfect reconstruction of dataset, thus proving the validity of the algorithm.

\begin{figure}[ht]
\includegraphics[width=0.9\linewidth,height=8cm]{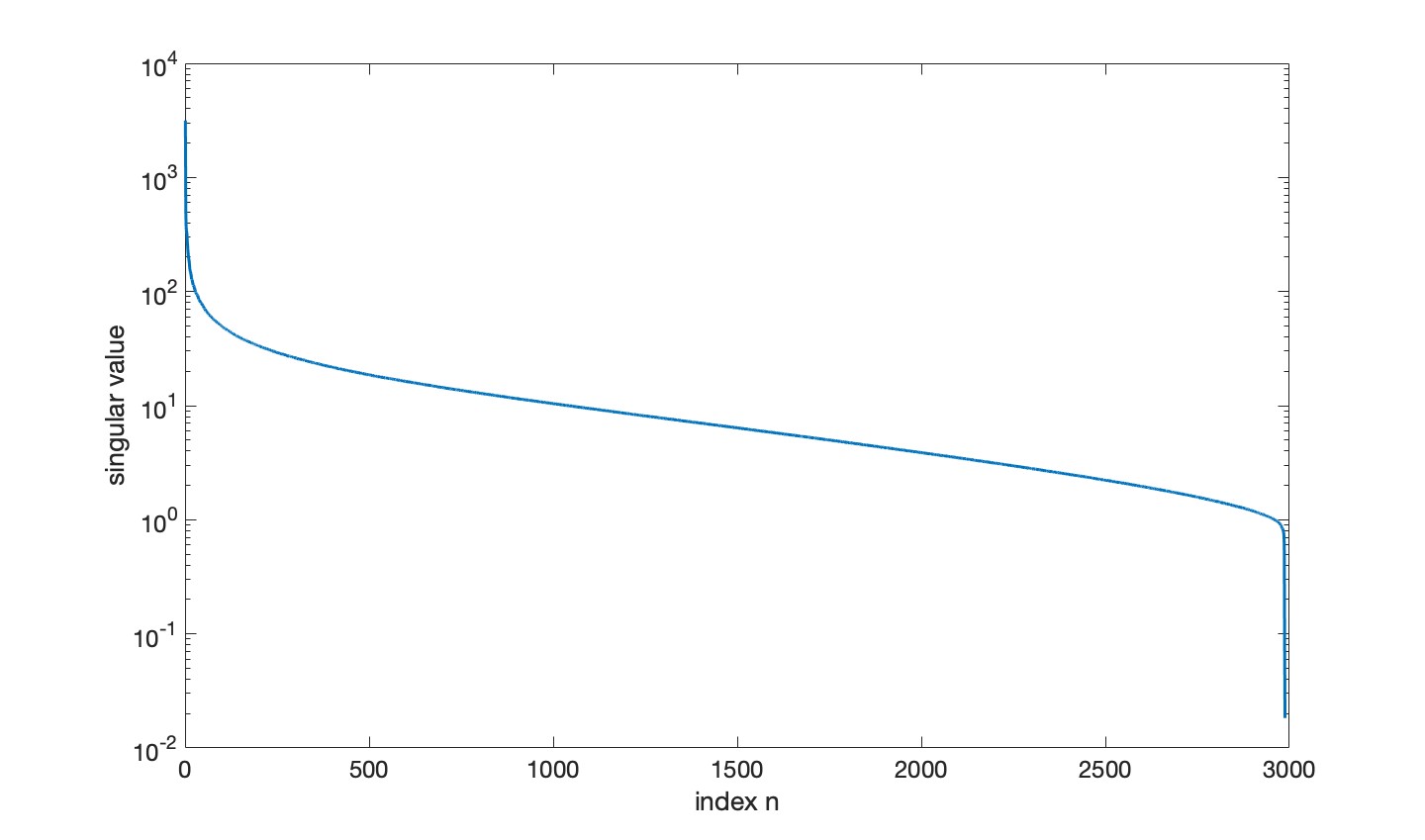}
\caption{\large{The behavior of singular values as obtained from CALTECH-101 dataset.}}
\label{fig:caltech101singularValues}
\end{figure}

\begin{figure}[ht]
\includegraphics[width=0.9\linewidth,height=8cm]{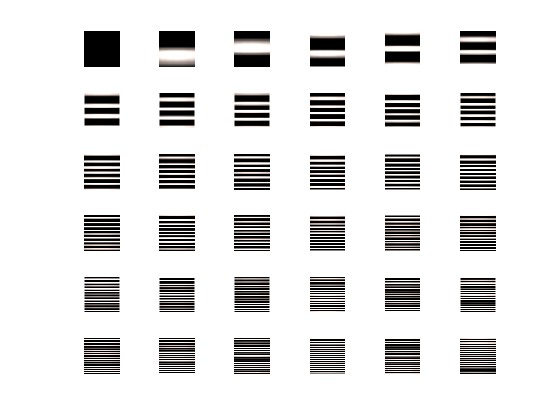}
\caption{\large{some basis elements from CALTECH-101 dataset corresponding to the first eigenvalues.}}
\label{fig:caltech101basis}
\end{figure}

Fig.~\ref{fig:caltech101singularValues} shows the behavior of singular values as obtained from CALTECH-101 dataset,
while some basis elements from CALTECH-101 dataset corresponding to the first singular values are reported in Fig.~\ref{fig:caltech101basis}.

A truncated representation of dataset has been achieved from (\ref{eq4B14}) by retaining only the components corresponding to the most significant singular values in matrix $S$. 

\begin{figure}[ht]
\includegraphics[width=0.9\linewidth,height=8cm]{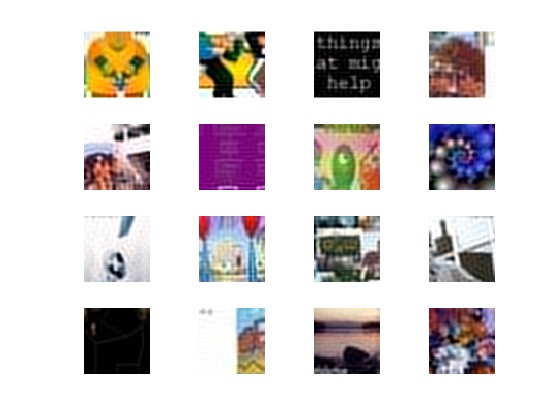}
\caption{\large{The same images of Fig.~\ref{fig:caltech101} reconstructed with truncated basis.}}
\label{fig:caltech101truncatedRecon}
\end{figure}

Fig.~\ref{fig:caltech101truncatedRecon} shows the same images as Fig.~\ref{fig:caltech101}, reconstructed with the 500 most significant components of the basis.

\subsection{Experiments on FLOWERS}

 The objective of the last experiment is to validate the approach reported in Subsection 4.C to derive a basis in a subspace. In this experiment the  FLOWERS dataset 
\cite{lim2005singular} that contains 4242 colour images of flowers has been used. 
The pictures are divided into five classes: chamomile, tulip, rose, sunflower, dandelion. 
For each class there are about 800 photos. Photos, which are not high resolution (about 320x240 pixels),  are  reduced to a single size
of ${128 \times 128 \times 3}$.
In this experiment  we use 3660 RGB images of the set for image reconstruction by the representation (\ref{eq4.1.1}),  with the basis achieved solving the eigenvalue problem discussed in Sect.4.C.
For this purpose, the data has been organized as a order-$4$ tensor $\mathtens{X} = 
\{ \mathtens{X}_{\textbf{i},n},n=1:N\}$, with $N = 3660, \; \mathtens{X} \in \mathbb{R}^{128 \times 128 \times 3\times3660}$, so that the dimension of the basis is $L = 128 \cdot 128 \cdot 3 = 49152$.

\begin{figure}[ht]
\includegraphics[width=0.9\linewidth,height=8cm]{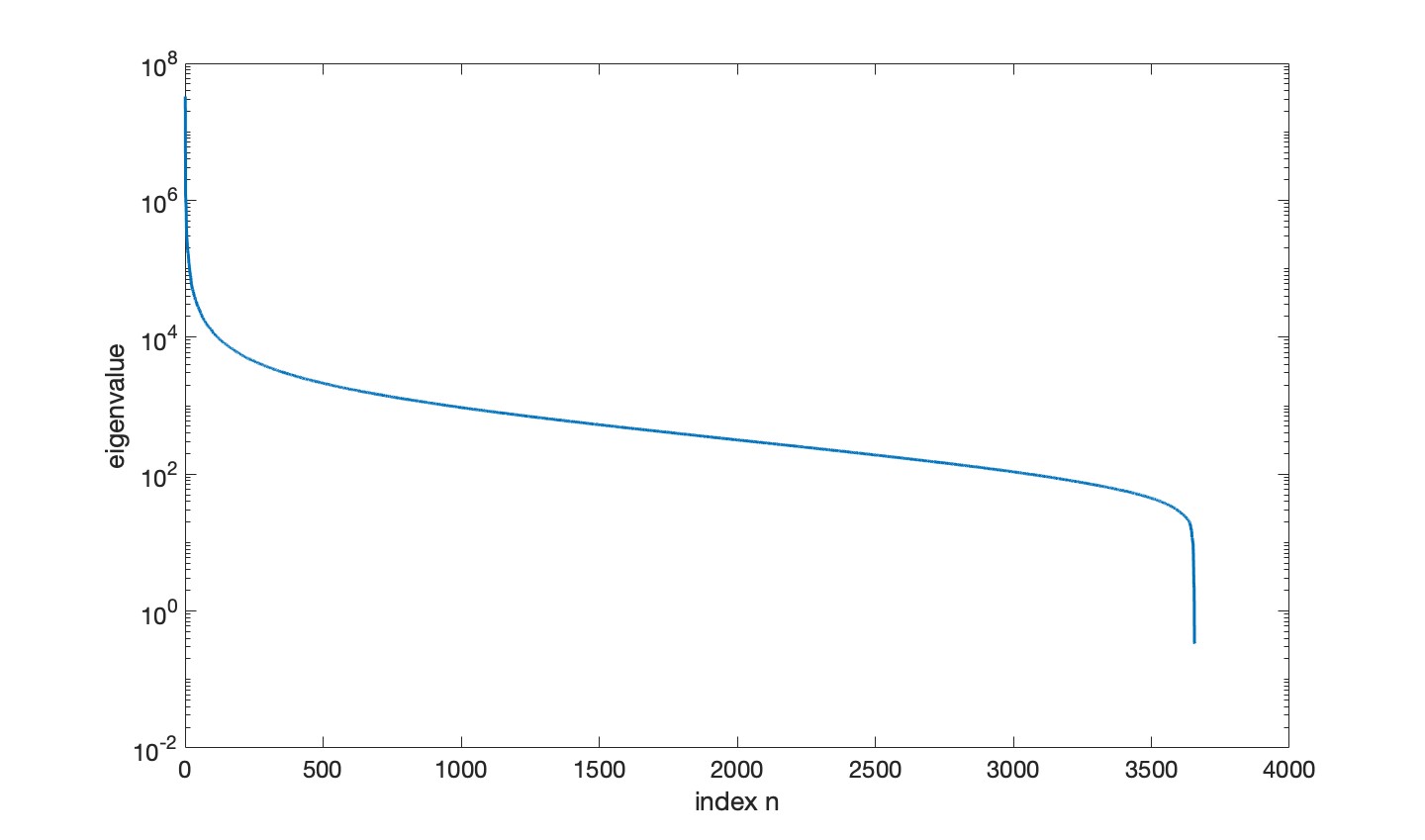}
\caption{\large{The behavior of eigenvalues as obtained from FLOWERS dataset.}}
\label{fig: FLOWERSeigenvalues}
\end{figure}

\begin{figure}[ht]
\includegraphics[width=0.9\linewidth,height=8cm]
{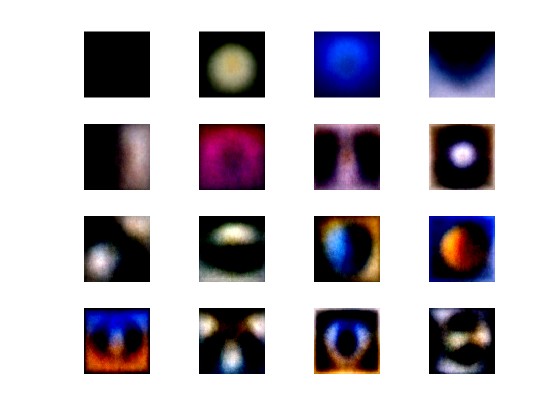}
\caption{\large{Some basis elements from FLOWERS dataset corresponding to the first eigenvalues.}}
\label{fig:FLOWERSbasis}
\end{figure}

\begin{figure}[ht]
\includegraphics[width=0.9\linewidth,height=8cm]
{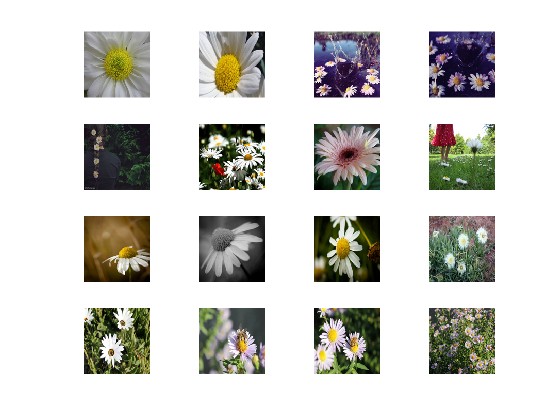}
\caption{\large{Several images taken from the dataset FLOWERS.}}
\label{fig:FLOWERS}
\end{figure}

\begin{figure}[ht]
\includegraphics[width=0.9\linewidth,height=8cm]{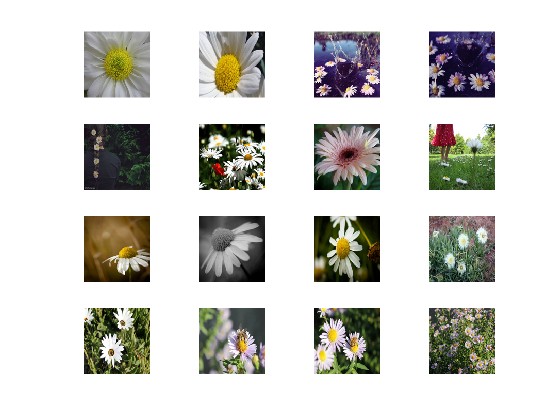}
\caption{\large{The same images of Fig.~\ref{fig:FLOWERS}reconstructed by (\ref{eq4.1.1}).}}
\label{fig:FLOWERSreconstructed}
\end{figure}

\begin{figure}[ht]
\includegraphics[width=0.9\linewidth,height=8cm]{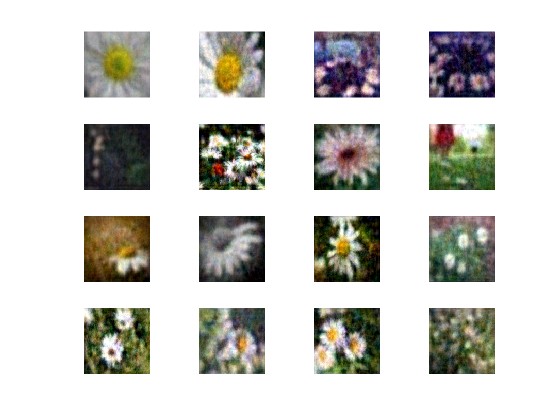}
\caption{\large{The same images of Fig.~\ref{fig:FLOWERS} reconstructed with truncated basis.}}
\label{fig:FLOWERStruncatedRecon}
\end{figure}

Fig.~\ref{fig: FLOWERSeigenvalues} shows the behavior of eigenvalues as obtained from FLOWERS dataset, using Algorithm 3,
while some basis elements from FLOWERS dataset corresponding to the first eigenvalues are reported in Fig.~\ref{fig:FLOWERSbasis}.
Fig.~\ref{fig:FLOWERS} reports several images taken from the dataset.  For comparison the same images reconstructed by (\ref{eq4.1.1}), with the basis achieved with Algorithm 3, are reported in Fig.~\ref{fig:FLOWERSreconstructed}, showing a perfect match. These results confirm the validity of the Algorithm 3 to derive a basis of subspace.
A truncated representation of dataset has been achieved from (\ref{eq4C14}), by retaining only the components corresponding to the most significant eigenvalues in matrix $\hat\Sigma$.
Fig.~\ref{fig:FLOWERStruncatedRecon} shows the same images as Fig.~\ref{fig:FLOWERS}, reconstructed with the 300 most significant components of the basis.

\section{Conclusion}
Tensor PCA is a multilinear extension of PCA, one of the most popular and efficient subspace techniques used in matrix analysis to reduce the dimensionality of a dataset. In this paper a mathematical framework for tensor PCA has been developed by deriving a basis in tensor space and by projecting data on a subspace. In particular the problem of deriving a basis from a self-adjoint tensor operator has been studied and the equivalence between eigenvalue equation for this operator and standard  matrix eigenvalue equation has been proven.

%

\newpage
\clearpage
%
%

\end{document}